\pdfoutput=1

\RequirePackage{fix-cm}

\documentclass[leqno, fleqn, centertags, 12pt]{article}   

\usepackage{latexsym}
\usepackage{amsmath}
\usepackage{amscd}
\usepackage{amssymb}
\usepackage{verbatim}

\usepackage[T1]{fontenc}
\usepackage{latexsym}
\usepackage{manfnt}
\usepackage[sans]{dsfont}
\usepackage{palatino}
\usepackage[sc, osf]{mathpazo} 
\usepackage{eulervm}

\usepackage{fullpage}

\usepackage{graphicx}

\newcommand{\nnn}{\mathbb{N}}
\newcommand{\zzz}{\mathbb{Z}}
\newcommand{\id}{\mathop{\mathrm{id}}\nolimits}
\newcommand{\card}{\mathop{\mathrm{card}}\nolimits}

\newcommand{\minus}{\hspace*{0.15em}\mbox{\rule[0.4ex]{0.4em}{0.4pt}}\hspace*{0.15em}}

\newcommand{\ccomp}{\mathsf{c}}

\def\sss{\hskip.05em\ }
\def\dss{\hskip.1em\ }
\def\trs{\hskip.15em\ }
\def\qss{\hskip.2em\ }
\def\oss{\hskip.4em\ }

\def\halfff{{\hskip.025em}}
\def\fff{{\hskip.05em}}
\def\dff{{\hskip.1em}}
\def\trf{{\hskip.15em}}
\def\qff{{\hskip.2em}}
\def\off{{\hskip.4em}}

\def\ffdot{\hspace*{-0.1em}.\hspace*{0.2em}\ }
\def\dfdot{\hspace*{-0.2em}.\hspace*{0.4em}\ }

\def\ffcom{\hspace*{-0.1em},\hspace*{0.2em}\ }
\def\dfcom{\hspace*{-0.2em},\hspace*{0.4em}\ }

\newcommand{\nsp}{\hspace*{-0.1em}}
\newcommand{\dnsp}{\hspace*{-0.2em}}

\parskip=\the\bigskipamount

\makeatletter
\renewcommand{\@makefntext}[1]{\vspace*{0.5ex}\parindent=0em
\hspace*{-0.4em}
\hbox to 0.4em{\hss\@makefnmark}\hspace*{0.4em}{#1}
}
\makeatother

\newcounter{mysectionnumber}
\newcommand{\mysection}[2]{\setcounter{footnote}{0}
\refstepcounter{mysectionnumber}
\section*{ \textnormal{{\themysectionnumber.} {#1}}}\label{#2}}

\newcommand{\mynonumsection}[1]{\setcounter{footnote}{0}
\section*{\textnormal{#1}}}

\newcounter{myparnum}[mysectionnumber]
\newcommand{\mypar}[2]{\refstepcounter{myparnum}{\vspace{\bigskipamount}\textbf{\textit{{\themyparnum. #1}\label{#2}}}\hspace{0.5em}}}
\renewcommand{\themyparnum}{\themysectionnumber.\arabic{myparnum}}

\newcounter{myapparnum}[mysectionnumber]
\newcommand{\myappar}[2]{\refstepcounter{myapparnum}{\vspace{\bigskipamount}\textbf{\textit{{\themyapparnum. #1}\label{#2}}}\hspace{0.5em}}}
\renewcommand{\themyapparnum}{A.\themyappendnumber.\arabic{myapparnum}}

\newcommand{\myitpar}[1]{\vspace{\bigskipamount}\textbf{\textit{#1}}\hspace*{0.5em}}
\newcommand{\myit}[1]{\textbf{\textit{#1}}\hspace{0.0em}} 
\newcommand{\mytitle}[1]{\textbf{\textit{#1}}}

\newcounter{mylemmanum}[myparnum]
\newcommand{\proof}{\vspace{\bigskipamount}{\textbf{{\emph{Proof}.}}\hspace*{0.7em}}}
\newcommand{\prooftitle}[1]{\vspace*{\bigskipamount}{\textbf{{\emph{#1}.}}\hspace{0.7em}}}
\newcommand{\eproof}{ $\blacksquare$}
\newcommand{\esubproof}{ $\square$}

\newcounter{myappendnumber}
\newcommand{\myappend}[2]{\setcounter{footnote}{0}
\setcounter{myapparnum}{0}
\refstepcounter{myappendnumber}
\section*{\textnormal{Appendix {\themyappendnumber.} {#1}}}\label{#2}}  

\newcounter{myaparnum}[myappendnumber]
\newcommand{\ppo}[1]{\mathcal{P}\dff({#1})}
\newcommand{\bbb}{\mathcal{B}}

\newcommand{\alm}{\mathcal{A}}
\newcommand{\ove}{\mathcal{O}}

\newcommand{\cbbb}{\mathcal{B^\ccomp}}
\newcommand{\csup}{{^\ccomp}}

\newcommand{\comm}{\dff,\qff}
\newcommand{\aazz}{a\dff,\qff z}
\newcommand{\aabb}{a\dff,\qff b}
\newcommand{\sdiff}{\bigtriangleup}

\begin{document}

\title{The\qss Tutte\qss expansion\qss revisited}
\date{}
\author{\textnormal{Nikolai\hspace*{0.1em} V.\hspace*{0.1em} Ivanov}}

\footnotetext{\hspace*{-0.5em}\copyright\ Nikolai V. Ivanov, 2016.\trs 
Neither the work reported in this paper, nor its preparation were supported by any 
governmental or non-governmental agency, foundation, or institution.}

\maketitle

\vspace*{6ex}

\myit{\hspace*{0em}\large Contents}\vspace*{1ex} \vspace*{\bigskipamount}\\ 
\myit{Introduction}\hspace*{0.5em}  \hspace*{0.5em} \vspace*{1ex}\\
\myit{\phantom{1}1.}\hspace*{0.5em} Pre-matroids and matroids  \hspace*{0.5em} \vspace*{0.25ex}\\
\myit{\phantom{1}2.}\hspace*{0.5em} Triangles  \hspace*{0.5em} \vspace*{0.25ex}\\
\myit{\phantom{1}3.}\hspace*{0.5em} Permutations and transpositions  \hspace*{0.5em} \vspace*{0.25ex}\\
\myit{\phantom{1}4.}\hspace*{0.5em} Linkings  \hspace*{0.5em} \vspace*{0.25ex}\\
\myit{\phantom{1}5.}\hspace*{0.5em} Orders  \hspace*{0.5em} \vspace*{0.25ex}\\
\myit{\phantom{1}6.}\hspace*{0.5em} Multi-sets  \hspace*{0.5em} \vspace*{0.25ex}\\
\myit{\phantom{1}7.}\hspace*{0.5em} The Tutte polynomial and the Whitney multiset  \hspace*{0.5em} \vspace*{0.25ex}\\
\myit{\phantom{1}8.}\hspace*{0.5em} Branching and balance  \hspace*{0.5em} \vspace*{0.25ex}\\  
\myit{\phantom{1}9.}\hspace*{0.5em} Forced balance  \hspace*{0.5em} \vspace*{0.25ex}\\
\myit{10.}\hspace*{0.5em} Coda: the order-independence  \hspace*{0.5em} \vspace*{1ex}\\
\myit{Appendix 1.}\hspace*{0.5em} The symmetric exchange property  \hspace*{0.5em} \vspace*{0.25ex}\\
\myit{Appendix 2.}\hspace*{0.5em} Duality  \hspace*{0.5em} \vspace*{0.25ex}\\
\myit{Appendix 3.}\hspace*{0.5em} Permutations and transpositions  \hspace*{0.5em} \vspace*{0.25ex}\\
\myit{Appendix 4.}\hspace*{0.5em} A direct proof of Lemma \ref{xy-u-sets}  \hspace*{0.5em} \hspace*{0.5em} \vspace*{0.25ex}\\
\myit{Appendix 5.}\hspace*{0.5em} Classification of linkings  \vspace*{1ex}\\
\myit{Note historique}\hspace*{0.5em}   \hspace*{0.5em}  \vspace*{0.25ex}\\
\myit{Bibliographie}\hspace*{0.5em}  \hspace*{0.5em}  \vspace*{0.25ex}

\vspace*{\bigskipamount}
{\small
With the exception of Introduction,\qss the present paper is self-contained modulo basic concepts related to sets and maps.\qss
In particular{},\qss no knowledge of the matroid theory or of the graph theory is assumed.\qss

}

\newpage

\renewcommand{\baselinestretch}{1.01}
\selectfont

\mynonumsection{Introduction}

\myitpar{The Tutte polynomial and the Tutte expansion.}\qss This paper is a result of a reflection 
upon a classical theorem of\qss W.T. Tutte \cite{t1}.\qss
Let $G$ be a connected graph.\qss
Suppose that a linear order on the set $E$ of edges of $G$ is given.\qss
W.T. Tutte \cite{t1} defined a polynomial $\chi\fff(G,\dff\mathbold{x},\dff\mathbold{y})$ in two variables 
$\mathbold{x},\dff\mathbold{y}$ as the sum\vspace*{\medskipamount} 
\begin{equation*}
\quad
\chi\fff(G;\qff\mathbold{x},\dff\mathbold{y})\qff\qff =\qff\qff 
\sum_T {\dff}\mathbold{x}^{r(T)}\dff \mathbold{y}^{s(T)},
\end{equation*}
where $T$ runs over the set $\mathcal{T}$ of all spanning trees of $G$\dfdot
The polynomial $\chi\fff(G;\qff \mathbold{x},\qff \mathbold{y})$ 
is called\qss \emph{the Tutte polynomial of}\sss $G$\dfcom
and the above sum
is called\qss \emph{the Tutte expansion of}\qss $\chi\fff(G;\qff\mathbold{x},\dff\mathbold{y})$\dnsp.\qss
The natural numbers $r(T)$ and $s(T)$ are the so-called\qss \emph{internal}\qss and\qss \emph{external activities}\dss
of a spanning tree $T$ \emph{with respect to the given linear order} on $E$\dfdot
The internal and external activities of a spanning tree depend on the choice of a linear order on $E$\dfcom
but the polynomial $\chi\fff(G;\qff\mathbold{x},\dff\mathbold{y})$ 
{does not depend on this choice}.\qss 
Therefore $\chi\fff(G;\qff\mathbold{x},\dff\mathbold{y})$ is an invariant of the graph $G$\dfdot
This independence on the choice of linear order is a striking result,\qss 
for which Tutte gave a beautiful and intriguing proof.\qss
We will call this result,\qss as also its generalization to matroids (see below),\qss
the\qss \emph{Tutte order-independence theorem}.\qss

\myitpar{The paper.} The present paper is devoted to an elementary{\nsp},\qss detailed,\qss 
and self-contained proof of the Tutte order-independence theorem.\qss
The remaining part of Introduction is devoted to a discussion of the motivation 
behind and the novel aspects of this proof.\qss

\myitpar{The maps $\varphi,\qff\qff\qff\psi$\nsp.}\qss The internal and external activities $r(T)$ and $s(T)$ 
of a spanning tree $T$ in a graph are defined as the numbers of\qss 
\emph{internally}\qss and,\qss respectively{\nsp},\qss \emph{externally active}\qss
edges of $G$ with respect to a given order on the set $E$ of edges of $G$\dfdot
The Tutte definition of internally and externally active edges appears to be rather idiosyncratic.\qss 
In the present paper we do not use these notions at all,\qss
and define the Tutte polynomial in terms of the following two maps 
having spanning trees as values.\qss

Let $\alm$ and $\ove$ be the sets of spanning subgraphs resulting\fff,\qss
respectively{\nsp},\qss from removing an edge from a spanning tree and adding an edge to a spanning tree.\qss
A linear order $<$ on the set of edges allows to define 
natural maps $\varphi\colon\alm\qff\longrightarrow\qff\mathcal{T}$ 
and $\psi\colon\ove\qff\longrightarrow\qff\mathcal{T}$\nsp\dfdot
Namely{\nsp},\qss 
\[\quad
\varphi(D)\qff\qff =\qff\qff D\qff +\qff x
\hspace*{0.8em}\mbox{ and }\hspace*{0.8em}
\psi(Q)\qff\qff =\qff\qff Q\qff -\qff y,
\]
where $x$ is the $<$\dnsp-minimal element of $E\dff\smallsetminus\dff D$ such that $D\qff +\qff x$ is a spanning tree\dfcom
and $y$ is the $<$\dnsp-minimal element of $Q$ such that $Q\qff -\qff y$ is a spanning tree\dfcom
and we denote by plus and minus the operations of adding an edge to a subgraph 
and removing an edge from a subgraph,\qss respectively{\nsp}.\qss
For a spanning tree $T$ the numbers of elements in preimages $\varphi^{-1}\dff(T)$ and $\psi^{-1}\dff(T)$
are equal to the internal and external activities of $T$\dfcom respectively{\nsp},\qss
as one can easily check.\qss
The maps
\[ 
\quad
\varphi\dff\colon\dff\alm\qff\longrightarrow\qff\mathcal{T}\hspace*{1em}
\mbox{ and }\hspace*{1em}\psi\dff\colon\dff\ove\qff\longrightarrow\qff\mathcal{T}
\] 
are the central characters in our approach to the Tutte order{-}independence theorem.\qff\qss

\myitpar{The symmetry of the proof.}\qss The heart of Tutte's proof of the order-independence theorem is an analysis of 
the effect of replacing a given linear order by a new one 
differing from it only by the order of two consecutive edges.\qss
The non-trivial part of this analysis splits into four cases requiring 
four similar{},\qss but independent,\qss arguments.\qss
These four cases exhibit a striking $\mu_2\times\mu_2$ symmetry{\nsp},\qss
where $\mu_2=\{\dff{\minus}1\dff,\qff 1\dff\}$ is the multiplicative group of square roots of $1$\dfdot 

It turns out that the basic topological tools of the theory of graphs,\qss
namely{\nsp},\qss connected components and cycles,\qss 
only obscure Tutte's proof and partially destroy its symmetry{\nsp}.\qss
Apparently{\nsp},\qss it is easier to deal with connected components than with cycles,\qss
and by this reason connected components are preferred to cycles even when considering cycles 
is more natural from a logical point of view.\qss
At the same time,\qss Tutte proof uses pictures in his proof,\qss
and,\qss apparently{\nsp},\qss prefers more intuitive arguments 
involving connected components to arguments involving cycles.\qss

\myitpar{Matroids.}\qss An attempt to at least partially understand this 
$\mu_2\times\mu_2$ symmetry within our framework 
inevitably leads to the realization that the Tutte polynomial 
and the Tutte order-independence theorem are not really about graphs.\qss
The right framework for Tutte's results is the theory of matroids.\qss

Let us define a\qss \emph{pre-matroid on a set}\dss $X$ as a non-empty collection of subsets of $X$\dfcom
called\qss \emph{bases},\qss such that no basis properly contains another basis.\qss
A pre-matroid is called a\qss \emph{matroid}\qss 
if it satisfies\qss \emph{the exchange property}\qss (see Section \ref{basic}).\qss
Let us denote by $E\dff(H)$ the set of edges of a graph $H$\dfdot
The matroid associated to the graph $G$ is a matroid on the set $E\qff =\qff E\dff(G)$ having
as bases all subsets of the form $E\dff(T)\qff\subset\qff E$\dfcom 
where $T\dff\in\dff\mathcal{T}$ is a spanning tree of $G$\dfdot
By the very definition,\qss the Tutte polynomial $\chi\fff(G;\qff\mathbold{x},\dff\mathbold{y})$
depends only on this matroid associated with $G$\dfcom
but Tutte's proof of his order-independence theorem depends also on the graph $G$ behind this matroid.\qss

We will deal with matroids from the very beginning\fff.
In particular{}, we will prove the Tutte order-independence theorem for matroids.\qss
As a side benefit,\qss our proof disentangles purely combinatorial ideas of Tutte 
from topological arguments of the theory of graphs.\qss
The latter are replaced by topology-independent properties of matroids.\qss

The definition of the Tutte polynomials was extended to matroids long ago by H. Crapo \cite{c},\qss
whose motivation,\qss apparently{\nsp},\qss was quite different.\qss
Namely{\nsp},\qss Crapo's paper is based on his Ph.D. thesis written under supervision of
G.-{\dff}C. Rota and appears to be a part of a far reaching program of G.-{\dff}C. Rota 
of integrating combinatorics into the so-called mainstream mathematics.\qss
The author hopes that this paper will serve the same goal.\qss

\myitpar{The symmetry of the proof and linked matroids.}\qss The concept 
of\qss \emph{matroids duality}\qss (see Appendix \ref{duality})\qss explains a half 
of the $\mu_2\times\mu_2$ symmetry of Tutte's proof,\qss
but leaves another half unexplained.\qss
In order to explain the whole $\mu_2\times\mu_2$ symmetry{\nsp},\qss we introduce notion of\qss 
\emph{linking}\qss between two matroids on the same set.\qss
While every matroid is linked only to itself and to its dual matroid,\qss
the notion of a linking identifies the essential features of the theory
and allows to replace all four similar arguments by a single one.\qss
The resulting theorem is Theorem\qss \ref{linked-almost-bases}\qss below,\qss 
the focal point of the present paper.\qss

Only the proof of Theorem\qss \ref{linked-almost-bases}\qss depends on the exchange property of matroids.\qss
All other arguments work for pre-matroids without any modifications
(and are usually stated as results about pre-matroids).\qss
The exchange property enters the proof of Theorem\qss \ref{linked-almost-bases}\qss only through 
Lemma\qss \ref{non-triangle}\qss and Lemma\qss \ref{triangle}.\qss
As a result,\qss Lemmas\qss \ref{non-triangle}\qss and\qss \ref{triangle}\qss emerge as
the combinatorial basis of the Tutte theory{\nsp}.\qss

\mysection{Pre-matroids\qss and\qss matroids}{basic}

\vspace*{\bigskipamount}
\emph{Let us fix once and for all a finite set $X${\halfff}\dnsp.}

\myitpar{Adding and deleting elements.} For a subset $Y\subset X$ 
and an element $x\dff\in\dff X$ not belonging to $Y$\dfcom
we denote by $Y\qff +\qff x$ the set $Y\cup\{\dff x\dff\}$\dfdot
Note that $Y\qff +\qff x$ is defined only if $x\dff\not\in\dff Y$\dfdot

Similarly, for a subset $Y\qff\subset\qff X$ and an element $y\qff\in\dff Y$ of this subset, 
we denote by $Y\qff -\qff y$ the set $Y\dff\smallsetminus\dff\{\dff y\dff\}$\dfdot
Note that $Y\qff -\qff y$ is defined only if $y\dff\in\dff Y$\dfdot

\myitpar{Pre-matroids.} Let $\ppo{X}$ be the set of all subsets of $X$\dfdot
A\qss \emph{pre-matroid},\qss or a\dss \emph{pre-matroid structure},\qss on $X$\dfcom
is simply a\qss \emph{non-empty}\qss subset\qss  $\bbb\qff\subset\qff\ppo{X}$\dfdot 
When a subset $\bbb\qff\subset\qff\ppo{X}$ is considered as a pre-matroid,\qss 
elements $B\dff\in\dff\bbb$ are called\qss \emph{bases}\qss of\qss $\bbb$\dfdot

\myitpar{Almost-bases.} Let $\bbb$ be a pre-matroid structure on $X$\dfdot
An\qss \emph{almost-basis}\qss of $\bbb$ is defined as a subset of $X$ 
of the form $B\qff -\qff x$\dfcom where $B\dff\in\dff\bbb$ and $x\dff\in\dff B$\dfdot
A subset $D\qff \subset\qff  X$ is an almost-basis if and only if 
$D\qff +\qff x$ is a basis for some element $x\dff\in\dff X\dff\smallsetminus\dff D$\dfdot
For an almost-basis $D$ we denote by $U\dff(D)$ the set of all  
$x\dff\in\dff X\dff\smallsetminus\dff D$ such that $D\qff +\qff x$ is a basis.\qss
In other terms,\qss $x\dff\in\dff U\dff(D)$ if and only if $x\dff\not\in\dff D$ and $D\qff +\qff x\in\bbb$\dfdot 
Suppose that $C\qff \subset\qff  X$\dfcom $x\comm\qff y\dff\not\in\dff C$\dfcom and $x\qff \neq\qff  y$\dfdot 
Then
\[
\quad
y\dff\in\dff U\dff(C\qff +\qff x)\hspace*{0.8em}\mbox{ \emph{if and only if} }\hspace*{0.8em}x\dff\in\dff U\dff(C\qff +\qff y).
\]
Indeed,\qss both $y\dff\in\dff U\dff(C\qff +\qff x)$ and $x\dff\in\dff U\dff(C\qff +\qff y)$ 
are equivalent to $C\qff +\qff x\qff +\qff y\dff\in\dff\bbb$\dfdot

\myitpar{Over-bases.} An\dss \emph{over-basis}\dss of $\bbb$ is defined as 
a subset of $X$ of the form $B\qff +\qff y$\dfcom where $B\dff\in\dff\bbb$ and $y\dff\not\in\dff B$\dfdot
A subset $Q\qff \subset\qff  X$ is an over-basis if and only if 
$Q\qff -\qff y$ is a basis for some element $y\dff\in\dff Q$\dfdot
For an over-basis $Q$ we denote by $C\dff(Q)$ the set of all  $y\dff\in\dff Q$ such that $Q\qff -\qff y$ is a basis.\qss
In other terms,\qss $x\dff\in\dff C\dff(Q)$ if and only if $x\dff\in\dff Q$ and $Q\qff -\qff x\dff\in\dff\bbb$\dfdot
Suppose that $Q\qff \subset\qff  X$\dfcom $x\comm y\dff\in\dff Q$\dfcom and $x\qff \neq\qff  y$\dfdot 
Then
\[\quad
y\dff\in\dff C\dff(Q\qff -\qff x)
\hspace*{0.8em}\mbox{ \emph{if and only if} }\hspace*{0.8em}
x\dff\in\dff C(Q\qff -\qff y).
\]
Indeed,\qss both $y\dff\in\dff C\dff(Q\qff -\qff x)$ and $x\dff\in\dff C\dff(Q\qff -\qff y)$ 
are equivalent to $Q\qff -\qff x\qff -\qff y\in\bbb$\dfdot

\myitpar{Matroids.} A pre-matroid $\bbb$ is called a\qss \emph{matroid}\qss 
{\qff}if the following\qss \emph{exchange property}\qss holds.\vspace*{-\bigskipamount}\vspace*{\medskipamount} 
\begin{quote}
\hspace*{1em}\emph{If\dff\qss $B_1\dff,\off B_2$\qss 
are bases of\dff\qss $\bbb$\qss and\qss $x\dff\in\dff B_1\dff\smallsetminus\dff B_2$\ffcom \\
\hspace*{0.8em}
then\qss $B_1\qff -\qff x\qff +\qff y$\qss is a basis for some\qss $y\dff\in\dff B_2\dff\smallsetminus\dff B_1$\nsp.}
\end{quote}
\vspace*{-\bigskipamount}\vspace*{\medskipamount}\nopagebreak[4] 
Recall that the\qss \emph{symmetric difference}\trs $P\sdiff Q$\qss of two sets\qss $P\comm Q$\qss is defined as
\[\quad
P\sdiff Q\off =\off (P\dff\smallsetminus\dff Q)\dff\cup\dff (Q\dff\smallsetminus\dff P).
\]
If $B\qff\subset\qff X$\ffcom $x\comm\qff y\dff\in\dff X$\ffcom and $x\dff\in\dff B$\dfcom $y\dff\not\in\dff B$\dfcom
then,\qss obviously,\qss $B\qff -\qff x\qff +\qff y\qff =\qff B\dff\sdiff\dff\{\dff x\comm y\dff\}$\dfdot

\mypar{Theorem (Symmetric exchange property).}{symm-exchange} \emph{If\qss $\bbb$\qss is a matroid,\dff\qss
$B_1\dff,\off B_2\in\bbb$\dfcom 
and\qss $x\in B_1\dff\smallsetminus\dff B_2$\ffcom 
then there exists\qss $y\in B_2\dff\smallsetminus\dff B_1$\qss 
such that\qss $B_1\sdiff\{\dff x\comm y\dff\}\fff\in\fff\bbb$\qss 
and\qss $B_2\sdiff\{\dff x\comm y\dff\}\fff\in\fff\bbb$\dnsp.}

\proof See Appendix \ref{symmetric exchange} for an elementary self-contained proof.  \eproof

\myitpar{Dual pre-matroids.} For a subset $Y\qff \subset\qff X$ we denote by $Y{\dff}\csup$ 
the complement $X\dff\smallsetminus\dff Y$ of $Y$ in $X$\dfdot
The\qss \emph{dual}\qss pre-matroid\qss $\bbb\csup$\qss of a pre-matroid\qss $\bbb$\qss is defined as
\[\quad
\bbb\csup\off = \off 
\bigl\{\dff B\csup \qff\mid\qff B\dff\in\dff\bbb \dff\bigr\}\fff\qff \subset\qff \ppo{X}.
\]
Obviously{\nsp},\qss $\bbb\csup\fff\csup\qff =\qff \bbb$\dfdot\newpage

\mypar{Theorem.}{dual-matroid-1} \emph{If\dss $\bbb$\dss is a matroid,\dss
then the dual pre-matroid\dss $\bbb\csup$\dss of\dss $\bbb$\dss is also a matroid.\dss}

\proof See Appendix \ref{duality}.\qss \eproof

\mysection{Triangles}{matroid-triangles}

\vspace*{\bigskipamount}
In this section we assume that $\bbb$ is a\qss \emph{matroid}\dff\qss structure on $X$\dfdot
The results of this section will be used only in Section \ref{linking-balance}.\qss

\mypar{Lemma.}{non-triangle} \emph{Suppose that\qss 
$a\comm z\dff\in\dff X$\qss and $a\qff \neq\qff  z$\ffcom 
and suppose the\qss $C\qff \subset\qff  X$\qss and\qss $a\comm z\dff\not\in\dff C$\dfdot 
Suppose that\qss $\bbb$\dss is matroid on\trs $X$\trs and\dff\qss 
$C\qff +\qff a\comm\qff\qff C\qff +\qff z$\qss are almost-bases of\qss $\bbb$\dfdot
In this situation either\dff\qss $C\qff +\qff a\qff +\qff z\dff\in\dff\bbb$\dfcom 
{\qff}or\qff\qss $U\dff(C\qff +\qff a)\qff =\qff U\dff(C\qff +\qff z)$\dfdot}

\proof Consider an arbitrary $e\in U\dff(C\qff +\qff a)$\dfdot
Then $C\qff +\qff a\qff +\qff e\dff\in\dff\bbb$\dfdot
Since $C\qff +\qff z$ is an almost-basis,\qss 
$C\qff +\qff z\qff +\qff v\dff\in\dff\bbb$ for some $v\dff\in\dff X$\dfdot
If $v\qff =\qff e$\dfcom then $e\qff =\qff v\dff\in\dff U\dff(C\qff +\qff z)$\dfdot
If $v\qff \neq\qff  e$\dfcom then 
$v\dff\in\dff (C\qff +\qff z\qff +\qff v)\dff\smallsetminus\dff (C\qff +\qff a\qff +\qff e)$
and the exchange property implies that
\[
C\qff +\qff z\qff +\qff b\off =\off (C\qff +\qff z\qff +\qff v)\qff -\qff v\qff +\qff b\dff\in\dff\bbb
\] 
for some $b\dff\in\dff (C\qff +\qff a\qff +\qff e)\dff\smallsetminus\dff (C\qff +\qff z\qff +\qff v)$\dfdot
Obviously{},\qss $b\qff =\qff a$\qss or\qss $e$\dfdot
If $b\qff =\qff a$\dfcom then $C\qff +\qff a\qff +\qff z\qff =\qff C\qff +\qff z\qff +\qff b\dff\in\dff\bbb$\dfdot
Otherwise,\qss $b\qff =\qff e$\dfcom and hence 
$C\qff +\qff z\qff + e\dff\in\dff\bbb$\dfcom i.e. $e\in U\dff(C\qff +\qff z)$\dfdot
Therefore,\qss if $C\qff +\qff a\qff +\qff z\dff\not\in\dff\bbb$\dfcom
then $e\dff\in\dff U\dff(C\qff +\qff z)$ for every $e\dff\in\dff U\dff(C\qff +\qff a)$\dfcom
i.e. $U\dff(C\qff +\qff a)\qff \subset\qff U\dff(C\qff +\qff z)$\dfdot
Similarly{},\qss if $C\qff +\qff a\qff +\qff z\dff\not\in\dff\bbb$\dfcom
then $e\dff\in\dff U\dff(C\qff +\qff a)$ for every $e\dff\in\dff U\dff(C\qff +\qff z)$\dfcom
i.e. $U\dff(C\qff +\qff z)\qff \subset\qff  U\dff(C\qff +\qff a)$\dfdot
The lemma follows.\qss  \eproof

\myitpar{Triangles.} A\qss \emph{triangle}\qss in $\bbb$ is 
a subset $C\qff \subset\qff  X$ together with three distinct elements\qss
$u\fff,\qff\qff v\fff,\qff\qff w\dff\in\dff X\dff\smallsetminus\dff C$\qss such that\qss 
$C\qff +\qff u\qff +\qff v\fff,\off 
C\qff +\qff u\qff +\qff w\fff,\off 
C\qff +\qff v\qff +\qff w\dff\in\dff\bbb$\dfdot

\mypar{Lemma.}{triangle} \emph{Suppose that\qss $C\subset X$\qss together with elements\qss
$a\fff,\qff\qff z\fff,\qff\qff d\dff\in\dff X\dff\smallsetminus\dff C$\qss is a triangle.\qss
If\qss $e\dff\in\dff U\dff(C+z)$\dfcom then either\qss 
$e\dff\in\dff U\dff(C+a)$\dfcom or\dff\qss $e\dff\in\dff U\dff(C+d)$\dfdot}

\proof If\qss $e\dff\in\dff U\dff(C\qff +\qff z)$\dfcom 
then $C\qff +\qff z\qff +\qff e\dff\in\dff\bbb$\dfdot 
Let us apply the exchange property to the bases 
$C\qff +\qff z\qff +\qff e$\dfcom $C\qff \qff +\qff a\qff +\qff d$ 
and 
\[
\hspace*{1em}
z\dff\in\dff (C\qff +\qff z\qff +\qff e)\dff\smallsetminus\dff (C\qff +\qff a\qff +\qff d).
\]
\newpage
Since $(C\qff +\qff a\qff +\qff d)\dff\smallsetminus\dff (C\qff +\qff z\qff +\qff e)\qff 
=\qff \{\dff a\dff,\off z\dff\}$\dfcom
the exchange property implies that either
\[
\hspace*{1em}C\qff +\qff a\qff +\qff e,
\hspace*{0.6em}\mbox{ or }\hspace*{0.6em}
C\qff +\qff d\qff +\qff e\dff\in\dff\bbb.
\]
In the first case $e\dff\in\dff U\dff(C+a)$\dfcom in the second case $e\dff\in\dff U(C+d)$\dfdot  \eproof

\mysection{Permutations\qss and\qss transpositions}{transpositions}

\myitpar{Permutations.} A \emph{permutation}\sss of $X$ is simply a bijection $X\to X$\dfdot
The set of all permutations of $X$ with the composition
$(\sigma\dff,\qff\tau\fff)\longmapsto\sigma\circ\tau$ as the binary operation is well known to be a group.\qss
By the very definition,\qss this group acts on $X$ on the left
and therefore acts on other sets canonically related to $X$\dfcom
for example on the set $\ppo{X}$ of all subsets of $X$\dfdot
Usually the composition $\sigma\circ\tau$ is denoted simply by $\sigma\dff\tau$\dfdot

\myitpar{Transpositions.} A permutation $\tau$ of $X$ is called\qss \emph{transposition}\qss 
if $\tau\qff \neq\qff \id_X$ and $\tau\dff(x)\qff =\qff x$ for all elements $x\dff\in\dff X$ except two.\qss 
Then $\tau$ interchanges these two elements.\qss
Clearly,\qss for any two distinct elements $\aabb$ of $X$
there is a unique transposition interchanging $a$ and $b$\ffdot
We will denote it by $\tau_{\fff a b}$\dfdot 
By the definition,\qss
\[
\quad
\tau_{\fff a b}(a)\off =\off b,\hspace*{1em} \tau_{\fff a b}(b)\off =\off a,
\hspace*{1em}\mbox{ and }\hspace*{1em}\tau_{\fff a b}(x)\off =\off x
\hspace*{0.7em}\mbox{ if }\hspace*{0.7em}x\qff \neq\qff  \aabb.
\]
The transposition $\tau_{\fff a b}$ is called the\qss \emph{transposition of}\qss $\aabb$\dfdot
Clearly{},\qss every transposition $\tau$ is a non-trivial involution,\oss 
i.e. $\tau\circ\tau\qff =\qff \id$ and $\tau\qff \neq\qff \id$\dfdot

If $\aabb\in X$\dfcom {\dnsp}$a\qff \neq\qff  b$\dfcom and $\sigma$ is a permutation of $X$\ffcom
{\qff}then $\sigma\circ\tau_{\fff a b}\circ\sigma^{-1}$ is the transposition 
of $\sigma(a)\fff,\off\sigma(b)$\dfcom as a trivial verification shows.\qss
Hence
\begin{equation}
\label{sigma-tau}
\quad
\tau_{\fff\sigma(a)\fff \sigma(b)}\off =\off  \sigma\circ\tau_{\fff a b}\circ\sigma^{-1}
\hspace*{1em}\mbox{ and }\hspace*{1em}
\tau_{\fff\sigma(a)\fff \sigma(b)}\circ\sigma\off =\off  \sigma\circ\tau_{\fff a b}.
\end{equation}

\vspace*{-\bigskipamount}
\myitpar{Action of transpositions on subsets of\qss $X$\nsp.}
Suppose that\qss $\aazz\dff\in\dff X$\qss and\qss $a\qff \neq\qff z$\ffcom 
and let\qss $\tau\qff =\qff \tau_{\fff a\dff z}$\ffdot
{\qff}Obviously,\oss $\tau\dff(Y)\qff =\qff Y$\qss for a subset $Y\subset X$ 
if and only if either both $\aazz\dff\in\dff Y$\dfcom or\dss both\qss $\aazz\dff\not\in\dff Y$\dfcom
and\qss $\tau\dff(Y)\qff \neq\qff Y$\qss if and only if exactly one of the elements $\aazz$\qss belongs to $Y$\dfdot
{\qff}If\qss $a\dff\in\dff Y$\qss and\qss $z\dff\not\in\dff Y$\dfcom then 
\[
\quad
\tau\dff(Y)\off =\off Y\qff -\qff a\qff +\qff z\fff,
\]
and if\qss $z\dff\in\dff Y$\qss and\qss $a\dff\not\in\dff Y$\dfcom then 
\[
\quad
\tau\dff(Y)\off =\off Y\qff -\qff z\qff +\qff a\fff.
\]

\newpage

\mysection{Linkings\qss of\qss pre-matroids}{linkings}

\myitpar{Linkings.} Let\qss $\bbb$\qss and\qss $\bbb^*$\qss be two pre-matroids on $X$\dfcom
and let\qss $L\dff\colon\bullet\qff\longmapsto\qff \bullet^*$\qss 
be a bijection\qss $\bbb\qff\longrightarrow\qff\bbb^*$\dfdot
The bijection\dss $L$\dss is said to be a\qss \emph{linking}\qss if\qss  
for every\qss $B\in\bbb$\qss and every transposition\qss $\tau\fff\colon X\qff \to\qff  X$\qss 
the following two\qss \emph{linking conditions}\qss hold.\vspace*{-\bigskipamount}\vspace*{\medskipamount}
\begin{quote}
$\mathbold{L1}$\dnsp.\fff\quad \emph{If\qff\qss $\tau\dff(B)\dff\in\dff\bbb$\dfcom
{\off}then\qff\qss $\tau\dff(B^*)\dff\in\dff\bbb^*$\qss {\qff}and\qff\qss $\tau\dff(B^*)\off =\off \tau\dff(B)^*$\dnsp.} 

$\mathbold{L2}$\nsp.\quad \emph{If\qff\qss $\tau\dff(B^*)\dff\in\dff\bbb^*$\dfcom
{\off}then\qff\qss $\tau\dff(B)\dff\in\dff\bbb$\qss {\qff}and\qff\qss $\tau\dff(B^*)\off =\off \tau\dff(B)^*$\dfdot} 
\end{quote} 

\vspace*{-\bigskipamount}\vspace*{\medskipamount}
Since\trs $L$\trs is a bijection,\qss
the condition $\mathbold{L2}$ is equivalent to the condition $\mathbold{L1}$ 
for the inverse map\qss $L^{\dff -1}$\dnsp.\dff\oss
In particular{},\qss a bijection\trs $L\dff\colon\bbb\qff\longrightarrow\qff\bbb^*$\trs 
is a linking if and only if\trs $L^{\dff -1}$\trs is a linking{\halfff}.\oss 
Obviously{\nsp},\qss the identity map\qss $\id\dff\colon\bbb\qff\longrightarrow\qff\bbb$\qss
and the map\qss $\ccomp\dff\colon\bbb\qff\longrightarrow\qff\bbb^*$\qss 
defined by\qss $\ccomp\dff\colon B\qff\longmapsto\qff B^\ccomp$\qss are linkings.\qss
The definition of a linking is motivated not by a desire of greater generality,\dff\qss
but by the desire to unify these two examples.\qss

\mypar{Theorem.}{linkings-classification}\qss \emph{If\dff\qss $L\dff\colon\bbb\qff \longrightarrow\qff \bbb^*$\qss 
is a linking{\halfff},\oss
then either\qss $\bbb^*\off =\off \bbb$\qss and\qss $L\off =\off \id$\nsp,\oss
or\oss $\bbb^*\off =\off \bbb^\ccomp$\oss and\oss $L\off =\off \ccomp$\ffdot}

\proof\qss See Appendix\qss \ref{app-linkings}.\qss  \eproof

\mypar{Lemma.}{xy-u-sets} \emph{Let\oss $L \dff\colon\dff \bbb\qff \longrightarrow\qff \bbb^*$\oss
be a linking{\halfff}.\oss
Suppose that\dss $S$\dss is an almost-basis of\oss $\bbb$\qss
and\dss $A$\dss is an almost-basis of\oss $\bbb^*$\dfdot
{\qff}Suppose that\dff\qss $x\fff,\off y\dff\in\dff U\dff(S)$\oss and\oss $x\qff \neq\qff y$\nsp.\linebreak
If\oss  $(S\qff +\qff x)^*\qff =\qff A\qff +\qff y$\ffcom
{\off}then\qss $x\dff\in\dff U^*(A)$\dnsp.}

\proof\qss By Theorem\qss \ref{linkings-classification},\oss it is sufficient
to consider only the linkings\qss $\id\dff\colon \bbb\qff \longrightarrow\qff  \bbb$\qss
and\qss $\ccomp\dff\colon \bbb\qff \longrightarrow\qff  \bbb\csup$\dnsp.\qff\oss
See Appendix\qss \ref{linkings-lemma-proof}\qss for a direct proof 
not relying on Theorem\qss \ref{linkings-classification}.\qss 

If\qss $\bbb^*\qff =\qff B$\qss and\qss $L\qff =\qff \id$\dnsp\halfff,\qff\oss
then\qss  $(S\qff +\qff x)^*\qff =\qff S\qff +\qff x$\nsp.\oss
But\qss $y\dff\not\in\dff S\qff +\qff x$\qss because\qss $y\qff \neq\qff x$\qss
and\qss $y\dff\not\in\dff S$\qss ({\fff}because\qss $S\dff \cap\dff U\dff (S)\qff =\qff \varnothing$\nsp)\halfff.\dff\oss 
Therefore,\qss in this case\qss $(S\qff +\qff x)^*$\qss cannot have the form\qss $A\qff +\qff y$\dfcom
and the lemma is trivially true.

If\qss $\bbb^*\qff =\qff \bbb\csup$\qss and\qss $L\qff =\qff\ccomp$\dnsp\halfff,\dff\oss
then\qss $A\qff +\qff y\qff =\qff (S\qff +\qff x)\csup$\qss and hence\qss
\[
A\off =\off (S\qff +\qff x)\csup\qff -\qff y\off 
=\off (S\qff +\qff x\qff +\qff y)\csup\fff,
\]
and\qss $A\qff +\qff x\qff =\qff (S\qff +\qff y)\csup\qff =\qff (S\qff +\qff y)^*$\dfdot
But\qss $S\qff +\qff y\fff\in\fff\bbb$\qss because $y\fff\in\fff U\dff(S)$\dfdot 
Therefore\qss
$A\qff +\qff x\qff =\qff (S\qff +\qff y)^*\dff\in\dff\bbb^*$\dfcom
and hence\qss $x\dff\in\dff U^*(A)$\dfdot  \eproof

\mysection{Orders}{orders}

\myitpar{Orders.} Given an order $\omega$ on $X$ and two elements 
$x\dff,\qff y\in X$\ffcom {\qff}both 
$\displaystyle x\trf  <_{\dff\omega} \qff y$
and $\displaystyle x\off _{\omega\dnsp}> \qff y$
will be used as a shorthand for\qss \emph{``{\dnsp}$x$\dss is less than\dss $y$\dss with respect to\dss $\omega$\dnsp''.}

We will consider only \emph{linear orders} on $X$\dfcom
i.e. orders $\omega$ on $X$ such that for every two elements $x\comm y\in X$ either $x=y$\dfcom
or $\displaystyle x\trf  <_{\dff\omega} \qff y$\dfcom
and $\displaystyle y\trf  <_{\dff\omega} \qff x$\dfdot
Given a linear order $\omega$ on $X$\dfcom
we will denote by $\min_{\dff\omega} Y$ the minimal element of a subset $Y\subset X$\dfdot

\myitpar{Action of permutations on orders.} Given an order $\omega$ on $X$ and a permutation $\omega\in\Sigma_{\dff X}$\ffcom
the order $\sigma\cdot\omega$ is defined as the unique order such that
\begin{equation}
\label{action}
\hspace*{1em}
\sigma\colon({\halfff} X\dff,\off {\omega} )\qff \to\qff ({\halfff} X\dff,\off {\sigma\cdot\omega} ) 
\end{equation}
is an isomorphism of ordered sets.\qss 
In other words,\qss if\dff\qss $x\fff,\qff y\in X$\dfcom then
\begin{equation*}
\hspace*{1em}
x\off  <_{\dff\omega} \off y\hspace*{1em}\mbox{ \emph{if and only if} }\hspace*{1em}\sigma\dff(x)\off <_{\dff\sigma\cdot\omega} \off \sigma\dff(y)\dff.
\end{equation*}
Obviously{\nsp}, $\sigma\cdot\omega$ is a linear order if and only if $\omega$ is.\qss
The map\qss $(\omega\fff,\trs\sigma)\longmapsto\sigma\cdot\omega$\qss is 
a left action of the group of permutation of $X$ on the set of orders on $X$\ffcom {\qff}i.e. 
\begin{equation}
\label{action-asso}
\hspace*{1em}
\tau\cdot(\sigma\cdot\omega)\off =\off (\tau\circ\sigma)\cdot\omega
\end{equation}
for all pairs\qss $\tau\fff,\qff\sigma$\qss of permutations of\qss $X$\qss and all orders $\omega$ on $X$\dfdot

\myitpar{Consecutive elements.} Two elements $x\dff,\dff y\dff\in\dff X$ are called 
\emph{consecutive elements with respect to the order} $\omega$
if either $x\qff  <_{\dff\omega} \qff y$ and there exist no elements $u\dff\in\dff X$ 
such that $x\qff <_{\dff\omega}  u \qff <_{\dff\omega} \qff y$\dfcom
or $y\qff <_{\dff\omega} \qff x$ and there exist no elements $u\dff\in\dff X$ 
such that $y\qff <_{\dff\omega} \qff u \qff <_{\dff\omega} \qff x$\dfdot

If $\omega$ is a linear order on $X$ and $\varepsilon$ is the transposition of two elements $\aazz$ 
consecutive with respect to $\omega$\ffcom
then the orders $\omega$ and $\varepsilon\cdot\omega$ differ only in order of elements $\aazz$\ffdot
In other words, $x\qff <_{\dff\omega} \qff y$ is equivalent to 
$x\qff <_{\dff\varepsilon\cdot\omega} \qff y$ unless $\{\dff x\comm y\dff\}\qff =\qff \{\dff\aazz\dff\}$\dfdot 
At the same time, $a\qff <_{\dff\omega} \qff z$ is equivalent to 
$z\qff <_{\dff\varepsilon\cdot\omega} \qff a$ by the definition of $\varepsilon\cdot\omega$\dfdot
It follows,\qss in particular,\qss 
that the elements $\aazz$ are consecutive with respect to the order $\varepsilon\cdot\omega$\dfdot

\myitpar{The graph of linear orders on $X$\dnsp.} The graph $\mathcal{L}_{{\halfff} X}$ of linear orders on $X$
has the set of all linear orders on $X$ as its set of vertices.\qss
Two linear order $\omega,\qff\omega'$ are connected by an edge if
$\omega'\qff =\qff \tau_{\fff a b}\cdot\omega$ for some elements $a\fff,\qff b\in X$ consecutive with respect to $\omega$\dfdot
If this is the case,\qss then also $\omega\qff =\qff \tau_{\fff a b}\cdot\omega'$ and
$a\fff,\qff b$ are consecutive with respect to $\omega'$\dfdot
Obviously{\nsp},\qss the transposition  $\tau_{\fff a b}$ is uniquely determined 
by the edge connecting $\omega,\qff\omega'$\nsp\dfcom
or{},\qss what is the same,\qss by the pair\qss  $\omega,\qff\omega'$\nsp\dfdot
For an edge $\mathcal{E}$ we will denote by $\tau\dff(\mathcal{E})$ the corresponding transposition.\qss

\myitpar{Remark.} The group of permutations of $X$ canonically acts on $\mathcal{L}_{{\halfff} X}$\dfdot
Indeed,\sss this group acts on the set of linear orders,\oss 
i.e. on the set of vertices of $\mathcal{L}_{{\halfff} X}$\ffdot
Let $\sigma$ be a permutation of $X$\ffdot
{\qff}Suppose that\qss $\omega'\qff =\qff \tau_{\fff a b}\cdot\omega$\qss 
for some elements $a\fff,\qff b\in X$ consecutive with respect to $\omega$\dfdot
By the second identity in\qss (\ref{sigma-tau})\qss 
\[
\quad
\sigma\cdot\omega'\off 
=\off \sigma\cdot\tau_{\fff\sigma(a)\fff \sigma(b)}\cdot\omega\off
=\off \tau_{\fff\sigma(a)\fff \sigma(b)}\cdot\sigma\cdot\omega.
\]
Obviously,\qss the elements $\sigma(a)\fff,\off \sigma(b)$ are consecutive with respect to $\sigma\cdot\omega$\dfdot
It follows that $\sigma$ takes edges to edges.\qss
Therefore,\qss the group of permutations of $X$ acts on $\mathcal{L}_{{\halfff} X}$\ffdot

\mypar{Lemma.}{connecting-pairs} \emph{Suppose that\qss $\omega\fff,\off \omega'$\qss are two linear orders on\qss $X$\dfdot
Then there exists a sequence
\[
\quad
\omega\qff =\qff \omega_{\fff 1}\fff,\hspace*{1.0em} \omega_{\fff 2}\fff,\hspace*{1.0em} \ldots\fff,\hspace*{1.0em} 
\omega_{\fff n}\qff =\qff \omega'
\]
of linear orders on\dss $X$\dss such that
for every\oss $i\off =\off 1\fff,\off 2\fff,\off \ldots\fff,\off n-1$\oss
the order\qss $\omega_{\fff i+1}$\qss is equal to\qss 
$\varepsilon_{\fff i}\cdot\omega_{\fff i}$\oss for some transposition\qss 
$\varepsilon_{\fff i}$\qss of two elements consecutive with respect to\qss $\omega_{\fff i}$\ffdot}

\proof\qss See Appendix \ref{perm-trans}.\qss  \eproof

\mypar{Corollary.}{connectedness} \emph{The  graph\qss $\mathcal{L}_X$\qss is connected.}\qss  \eproof

\mysection{Multi-sets}{multi-sets}

\myitpar{Multi-subsets.} A\qss \emph{multi-subset}\qss of a set $S$ 
is defined as a function $m\dff\colon\dff S\qff \longrightarrow\qff \nnn$\dfdot
The set $\{\dff s\dff\in\dff S\dff\mid\dff m(s)\qff \neq\qff 0\dff\}$ 
is called the\qss \emph{support}\qss of the multi-subset.\qss
A subset\qss $Y$\qss of\qss $S$\qss canonically defines a multi-subset of\qss $S$\dfcom
namely{\nsp},\qss it characteristic function\qss $\chi_Y\dff\colon\dff S\qff \longrightarrow\qff \nnn$\dfdot
Recall that $\chi_Y\dff(s)\qff =\qff 1$ if $s\dff\in\dff Y$ 
and $\chi_Y\dff(s)\qff =\qff 0$ if $s\dff\not\in\dff Y$\dfdot
Obviously{\nsp},\qss the support of the multi-set\qss $\chi_Y$\qss 
is nothing else but the subset $Y$\dfdot
By identifying subsets with their characteristic functions,\qss
we may consider subsets as multi-subsets.\qss

The values $m(s)$ of the characteristic function of a multi-subset 
are interpreted as multiplicities of element $s\dff\in\dff S$ in this multi-subset.\qss
In other words,\qss we think that a multi-subset $m$ 
contains $m(s)$ copies of $s$ for each $s\dff\in\dff S$\dfdot
Elements $s$ with $m(s)\qff =\qff 0$ are treated as 
elements not contained in the multi-subset\qss $m$\dfdot

\myitpar{Multi-sets.} A\qss \emph{multi-set}\qss is defined as a multi-subset 
of a fixed once and for all universal set $\mathcal{U}$\dfdot
Given a set\qss $S$\dfcom one can identify the multi-subsets of\qss $S$\qss
with multi-sets\qss $m$\qss
such that\qss $m(x)\qff =\qff 0$ for all $x\dff\not\in\dff S$\dfdot

\myitpar{Multi-images.} Let $f\dff\colon\dff S\qff \longrightarrow\qff  R$ be a map 
and $m\dff\colon\dff S\qff \to\qff \nnn$ be a multi-subset of $S$\dfdot
The multi-subset\qss 
$f\dff[\dff m\dff]\dff\colon\dff R\qff \longrightarrow\qff \nnn$\qss defined by\vspace*{0.5\medskipamount} 
\[
\quad
f\dff[\dff m\dff]\qff\colon\qff r\off\longmapsto\off 
\sum_{\qff f(s)\qff =\qff r}\dff m(s)\fff.
\]

\vspace{-1.5\medskipamount}
is called the\qss \emph{multi-image}\qss of $m$ under the map $f$\dfdot
For a subset $Y\qff \subset\qff  S$\dfcom 
the\qss \emph{multi-image}\qss $f\dff[\dff Y\dff]$\qss of\qss $Y$\qss
is defined as the multi-image $f\dff[\dff \chi_Y\dff]$\qss of the corresponding multi-subset.\qss 
Clearly{\nsp},\qss the multiplicity\qss of an element\qss $r\dff\in\dff R$ 
in the multi-image\qss $f\dff[\dff Y\dff]$\qss is equal to 
the number of elements of\qss $(\fff f\nsp\mid_{\dff Y\halfff})^{\dff -1}(r)$\dfcom
where\qss $f\nsp\mid_{\dff Y}\dff\colon Y\dff \to\qff R$\qss 
is the restriction of\qss $f$\qss to\qss $Y$\dfdot
The following lemma is obvious.\qss

\mypar{Lemma.}{multi-image-composition} \emph{If\qss $f\dff\colon\dff S\qff \to\qff  R$\dss 
is a map and\qss $g\dff\colon\dff S\qff \to\qff  S$\qss is a bijection,\qss
then} 
\[
\quad
f\dff[\dff S\dff]\off =\off f\circ g\dff[\dff S\dff].\hspace*{0.8em}\mbox{\eproof }
\]

\vspace*{-\bigskipamount}

\myitpar{Multi-subsets of\qff\qss $\nnn\times\nnn$\dnsp.}\oss 
Let\qss $\mathbold{x}\comm\mathbold{y}$\qss be two different variables.
One can associate with each multi-subset\qss $m\dff\colon \nnn\times\nnn\qff \longrightarrow\qff \nnn$ \qss
of\qss $\nnn\times\nnn$\qss a formal power series\qss $P_m(\mathbold{x}\comm\mathbold{y})$\qss 
in the variables\qss $\mathbold{x}\comm\mathbold{y}$\qss by the formula
\[
\quad
P_m\dff(\mathbold{x}\comm\mathbold{y})\off 
=\off \sum_{(a\comm b)\qff\in\qff\nnn\times\nnn} {\dff}m(a\comm b)\off \mathbold{x}^{a}\dff \mathbold{y}^{b}.
\]
Obviously{},\qss the map\qss $m\qff \longmapsto\qff P_m\dff(\mathbold{x}\comm\mathbold{y})$\qss establishes a one-to-one
correspondence between the multi-subsets of\qss $\nnn\times\nnn$\qss and formal power series with integer non-negative
coefficients in the variables\qss $\mathbold{x}\comm\mathbold{y}$\dfdot
Moreover,\qss the formal power series\qss $P_m(\mathbold{x}\comm\mathbold{y})$\qss is a polynomial 
if and only if the multi-subset\qss $m$\qss has finite support.\oss
The following lemma follows directly from the definitions.\qss

\mypar{Lemma.}{multi-image-in-nn}
\emph{Let\qss $S$\qss be a finite set and\qss $f\dff\colon\dff S\qff\to\qff \nnn\times\nnn$\qss be a map.\oss
Let\qss $a\comm b\dff\colon\dff S\qff\longrightarrow\qff\nnn$\qss be the two components of\qss $f$\ffcom
{\qff}i.e.\qss $f(s)\qff =\qff (a(s)\comm b(s))$\dff\qss for every\qss $s\dff\in\dff S$\dfdot
{\qff}Then the polynomial corresponding to the multi-image\qss $f\dff[\dff S\dff]$\qss is equal to}\vspace*{\medskipamount}
\begin{equation*}
\quad
P_{f\dff[\dff S\dff]}\dff(\mathbold{x}\comm\mathbold{y})\off 
=\off \sum_{s\dff\in\dff S}\off \mathbold{x}^{\dff a\fff(s)}\qff \mathbold{y}^{b\fff(s)}.\qquad \mbox{ \eproof }
\end{equation*}

\mysection{The\qss Tutte\qss polynomials\qss and\qss the\qss Whitney\qss multi-sets}{orders-pre-matroids}

\vspace*{\bigskipamount}
Let $\bbb$ be a pre-matroid on $X$
and let $\omega$ be a linear order on $X$\dfdot
Let $\alm$ and $\ove$ be the sets of all almost-bases and over-bases of $\bbb$\dfcom
respectively.\qss

\myitpar{The Tutte polynomial of a pre-matroid.}\qss Let\qss  
$\displaystyle
u_{\dff\omega}\dff\colon\dff\alm\qff\longrightarrow\qff X, 
\quad
c_{\fff\omega}\dff\colon\dff\ove\qff\longrightarrow\qff X
$\qss 
be the the maps defined,\qss respectively{\nsp},\dff\qss by 
\[
\quad
u_{\dff\omega}\dff(D)\off =\off \min\nolimits_{\dff\omega}\qff U\dff(D),\qquad\hspace*{0.6em} 
c_{\fff\omega}\dff(Q)\off =\off \min\nolimits_{\dff\omega}\qff C\dff(Q),
\]
Let\qss $
\varphi_{\omega}\dff\colon\dff\alm\qff\longrightarrow\qff\bbb,
\quad 
\psi_{\omega}\dff\colon\dff\ove\qff\longrightarrow\qff\bbb
$\qss 
be the the maps defined,\qss respectively{\nsp},\dff\qss by
\[
\quad
\varphi_{\omega}\dff(D)\off =\off D\off +\off u_{\dff\omega}\dff(D),\hspace*{2.0em}
\psi_{\omega}\dff(Q)\off =\off Q\off -\off c_{\fff\omega}\dff(Q).
\]
For a basis\qss $B\in\bbb$\dfcom {\qff}let\qss $i_{\dff\omega}\fff(B)$\qss and\qss $e_{\dff\omega}\fff(B)$ \qss
be the numbers of elements in the pre-images\qss
$(\varphi_{\omega})^{{\minus}1}\dff(B)$\qss and\qss $(\psi_{\omega})^{{\minus}1}\dff(B)$\qss respectively{\nsp}.\qss

Let\qss $\mathbold{x}\comm\qff\mathbold{y}$\qss be two different variables.\qss
The\qss \emph{Tutte polynomial}\qss of the pre-matroid\qss $\bbb$\qss 
with respect to the order\qss $\omega$\qss is defined as
\[
\quad
{T}_{\dff\omega}\dff(\fff\bbb\fff)\qff(\mathbold{x},\qff\mathbold{y})\off 
=\off \sum\nolimits_{\qff B\qff\in\qff\bbb}\off
\mathbold{x}^{\dff i_{\qff\omega}\fff(B)}\qff \mathbold{y}^{\fff e_{\dff\omega}\fff(B)}.
\]
It turns out that the Tutte polynomial\qss
$\displaystyle
{T}_{\dff\omega}\dff(\fff\bbb\fff)\qff(\mathbold{x},\qff\mathbold{y})
$\qss 
of\qss $\bbb$\qss does not depend on the order\qss $\omega$\qss if\qss $\bbb$\qss is a matroid.

\myitpar{Tutte's activities.} For a basis\qss $B\in\bbb$\dfcom {\qff}let\qss
$\displaystyle\mathcal{I}_{\omega}\fff(B)$\qss and\qss
$\displaystyle\mathcal{E}_{\omega}\fff(B)$ \qss be the sets
\[
\quad
\mathcal{I}_{\omega}\dff(B)\off =\off
\{\qff u_{\dff\omega}\dff(D) \qff\mid\qff D\dff\in\dff(\varphi_{\omega})^{{\minus}1}\dff(B)\qff\}\fff,\quad
\]
\[
\quad
\mathcal{E}_{\omega}\dff(B)\off =\off
\{\qff c_{\fff\omega}\dff(Q) \qff\mid\qff Q\dff\in\dff(\psi_{\omega})^{{\minus}1}\dff(B)\qff\}\fff.
\]

\vspace*{-\medskipamount}
By the very definition,\qss the sets\qss $\displaystyle\mathcal{I}_{\omega}\fff(B)$\qss and\qss
$\displaystyle\mathcal{E}_{\omega}\fff(B)$\qss are nothing else
but,\qss respectively,\qss the sets of\qss \emph{internally active}\qss 
and\qss \emph{externally active}\qss elements of\qss $X$\qss with respect
to the basis\qss $B\qss$ and the linear order\qss $\omega$\qss in the sense of Tutte.\oss
The maps\qss $u_{\dff\omega}$\qss and\qss $c_{\fff\omega}$\qss induce bijections
\[
\quad
(\varphi_{\omega})^{{\minus}1}\dff(B)\off \longrightarrow\off \mathcal{I}_{\omega}\dff(B),
\qquad
(\psi_{\omega})^{{\minus}1}\dff(B)\off \longrightarrow\off \mathcal{E}_{\omega}\dff(B),
\] 
and hence $i_{\dff\omega}\dff(B)$ and $e_{\dff\omega}\dff(B)$ are nothing else but,\qss
respectively,\qss the\qss \emph{internal activity}\qss and the\qss \emph{external activity}\qss 
of the basis\qss $B$\qss with respect to the order $\omega$ in the sense of Tutte.\qss
Therefore,\qss the polynomial\qss 
$\displaystyle {T}_{\dff\omega}\dff(\fff\bbb\fff)$\qss 
is nothing else but the classical Tutte polynomial.

\myitpar{The Tutte polynomial of a linking.} Let\qss $\bbb^*$\qss be another
pre-matroid on\qss $X$\qss and let\qss $\bullet\qff\longmapsto\qff\bullet^*$\dff\qss 
be a linking\qss $\bbb\qff\longrightarrow\qff\bbb^*$\dfdot
{\qff}Let\qss $\alm^*$\qss be the set of all almost-bases of\qss $\bbb^*$\dfdot

As before,\qss let\qss $\mathbold{x}\comm\qff\mathbold{y}$\qss be two different variables.\qss
The\qss \emph{Tutte polynomial}\qss of the linking\qss $\bbb\qff\longrightarrow\qff\bbb^*$\qss 
with respect to the order\qss $\omega$\qss is defined as
\[
\quad
{T}_{\dff\omega}\dff({\dff\bbb\qff\longrightarrow\qff\bbb^*})\qff(\mathbold{x},\qff\mathbold{y})\off 
=\off \sum\nolimits_{\qff B\qff\in\qff\bbb}\off
\mathbold{x}^{i_{\qff\omega}(B)}\qff \mathbold{y}^{i_{\dff\omega}(B^*)}.
\]

\vspace*{-\bigskipamount}
\myitpar{The Whitney multi-set of a linking.} The Tutte polynomial\qss
${T}_{\dff\omega}\dff({\dff\bbb\qff\to\qff\bbb^*})\qff(\mathbold{x},\qff\mathbold{y})$\qss
has a natural interpretation as a multi-subset of\qss $\nnn\times\nnn$\dfdot
Consider the map\qss 
\[ 
\quad
\mathfrak{W}_{\dff\omega}\dff({\dff\bbb\qff\longrightarrow\qff\bbb^*})\qff\colon\qff 
\bbb\qff\longrightarrow\qff\nnn\times\nnn
\] 
defined by\qss 
$\displaystyle
B\off\longmapsto\off \bigl(\dff i_{\dff\omega}\dff(B)\comm\off i_{\dff\omega}\dff(B^*) \dff\bigr)$\dfdot
The multi-image of $\bbb$ under this map is a multi-subset of\qss $\nnn\times\nnn$\dfdot
We call this multi-image the\qss \emph{Whitney multi-set}\qss of the linking\qss $\bbb\qff\to\qff\bbb^*$\qss 
with respect to the order\qss $\omega$\qss and denote it by\qss 
$\mathbb{W}_{\dff\omega}\qff({\dff\bbb\qff\to\qff\bbb^*})$\dfdot

\mypar{Theorem (Order-independence of the\dss Whitney multi-sets).}{w-independence}
\emph{If\dss $\bbb$\dss is a matroid, 
then the Whitney multi-set\qss 
$\displaystyle \mathbb{W}_{\dff\omega}\dff({\dff\bbb\qff\to\qff\bbb^*})\qff(\mathbold{x},\qff\mathbold{y})$\qss 
does not depend on the order\dss $\omega$\dfdot}\vspace*{\bigskipamount}

Sections\qss \ref{edges} -- \ref{proof}\qss are devoted to a proof of this theorem.

\mypar{Corollary\qss (Order-independence for linkings of matroids).}{main-theorem-linkings} 
\emph{If\qss $\bbb$\qss and $\bbb^*$ are matroids,\qss 
then the Tutte polynomial\qss 
$\displaystyle {T}_{\dff\omega}\dff({\dff\bbb\qff\to\qff\bbb^*})\qff(\mathbold{x},\qff\mathbold{y})$\qss 
does not depend on the order\qss $\omega$\dfdot}

\proof\qss Since\qss $\bbb\qff \subset\qff \ppo{X}$\qss is finite together with\qss $X$\dfcom
the multi-image\qss 
$\displaystyle \mathbb{W}_{\dff\omega}\qff({\dff\bbb\qff\to\qff\bbb^*})$\qss
of the map\qss 
$\displaystyle \mathfrak{W}_{\dff\omega}\dff({\dff\bbb\qff\to\qff\bbb^*})$\qss
has finite support{\halfff},\dff\qss
and hence the corresponding power series is a polynomial.\dff\qss
Lemma \ref{multi-image-in-nn} implies that this polynomial is equal to the Tutte polynomial
${T}_{\dff\omega}\dff({\dff\bbb\qff\to\qff\bbb^*})\qff(\mathbold{x},\qff\mathbold{y})$\dfdot
{\qff}It remains to apply Theorem \ref{w-independence}.\qss  \eproof

\mypar{Corollary\qss (Order-independence for matroids).}{main-theorem} \emph{If\qss $\bbb$\qss is a matroid,\oss 
then the Tutte polynomial\qss 
$\displaystyle {T}_{\dff\omega}\dff(\fff\bbb\fff)\qff(\mathbold{x},\qff\mathbold{y})$\qss 
does not depend on the order\dss $\omega$\dfdot}

\proof\qss Recall that\qss $\ccomp\colon B\qff \longmapsto\qff B\csup$\qss
is a linking\qss $\bbb\qff \longrightarrow\qff \bbb\csup$\dfdot
An immediate application of the matroid duality{\nsp}\qss 
(see Appendix \ref{duality},\qss Lemma \ref{phi-psi-dual}\fff)\qss shows that\qss
${T}_{\dff\omega}\dff({\dff\bbb\qff\to\qff\bbb\csup})\qff(\mathbold{x},\qff\mathbold{y})$\qss 
is equal to the Tutte polynomial\qss ${T}_{\dff\omega}\dff(\fff\bbb\fff)\qff(\mathbold{x},\qff\mathbold{y})$\dfdot
{\qff}It remains to apply Corollary \ref{main-theorem}.\qss  \eproof

\mysection{Branching\qss and\qss balance}{edges}

\myitpar{The framework.}\qss Let\qss $\bbb$\qss be pre-matroid on\qss $X$\dfcom
and let\qss $\alm$\qss be  the set of almost-bases of $\bbb$\ffdot
{\qff}Let\qss $\mathcal{E}$\qss be an edge of\qss $\mathcal{L}_{{\halfff} X}$\ffcom
and let\qss $\varepsilon\qff =\qff\tau\dff(\mathcal{E})$\qss be the corresponding transposition.\oss

Let\qss $\omega\fff,\off\pi$\qss be the linear orders on\qss $X$\qss connected by\qss $\mathcal{E}$\dfcom
{\qff}and let\qss $\aazz$\qss be the elements of\qss $X$\qss interchanged by\qss $\varepsilon$\dfdot
Then the orders\qss $\omega$\qss and\qss $\pi$\qss differ only by the order of elements\qss $\aazz$\ffcom
{\qff}and the elements\qss $\aazz$\qss are consecutive with respect to both\qss $\omega$\qss and\qss $\pi$\ffdot
{\qff}Without any loss of generality one may assume that\qss 
$a\dff  <_{\dff\omega}\dff z$\qss and\qss $z\dff  <_{\dff\pi\fff}\dff a$\halfff\dfdot

If\qss $A\dff\subset\dff X$\qss and\qss $\aazz\dff\not\in\dff A$\ffcom {\qff}then\qss
$\displaystyle \varepsilon(A\qff +\qff a)\off =\off A\qff +\qff z$\qss and\qss
$\displaystyle \varepsilon(A\qff +\qff z)\off =\off A\qff +\qff a$\dfdot
{\qff}In particular{},\qss
$\displaystyle \varepsilon(A\qff +\qff a)\off \neq\off A\qff +\qff a$\qss and\qss
$\displaystyle \varepsilon(A\qff +\qff z)\off \neq\off A\qff +\qff z$\dfdot

\myitpar{Branching almost-bases.}
An almost-basis\qss $A\dff\in\dff\alm$\qss is said to be\qss 
\emph{$\mathcal{E}$\dnsp-branching}\trf\qss if 
\begin{equation*}
\quad
\varphi_\omega(A)\off \neq\off \varphi_{\fff\pi\fff}(A).
\end{equation*}
Clearly{\nsp},\qss if $A$ is $\mathcal{E}$\dnsp-branching{\halfff},\qss
then the orders $\omega$ and $\pi$ differ on $U\dff(A)$ and hence $\aazz\dff\in\dff U\dff(A)$\dfdot
In particular{},\qss $\aazz\dff\not\in\dff A$\qss and hence\qss $\varepsilon\fff(A)\qff =\qff A$\dfdot
Moreover{},\qss one of the elements $\aazz$ should be equal to\qss $\min\nolimits_{\dff\omega} U\fff(A)$\dfcom
and the other to\qss $\min\nolimits_{\qff\pi} U\fff(A)$\dfdot
Since\qss $a\dff  <_{\dff\omega}\dff z$\qss and\qss $z\dff  <_{\dff\pi\fff}\dff a$\ffcom in this case
\begin{equation}
\label{branching-phi-image}
\quad
\min\nolimits_{\dff\omega} U\fff(A)\off =\off a
\hspace*{0.7em}\mbox{ {and} }\hspace*{0.7em}
\varphi_\omega(A)\off\halfff =\off A\qff +\qff a\fff,
\end{equation}
\begin{equation}
\label{branching-pi-image}
\quad
\min\nolimits_{\qff\pi\dff} U\fff(A)\off\halfff =\off\trf z
\hspace*{0.7em}\mbox{ {and} }\hspace*{0.7em}
\varphi_{\dff\pi\fff}(A)\off =\off A\qff +\qff\trf z\fff.
\end{equation}

\vspace*{-0.7\medskipamount}
Conversely{\nsp},\oss suppose that\qss $\aazz$\qss are 
the two smallest elements of\qss $U\fff(A)$\qss with respect to\qss $\omega$\dfdot
Since\qss $a\dff  <_{\dff\omega}\dff z$\qss and\qss $z\dff  <_{\dff\pi\fff}\dff a$\ffcom in this case
$\displaystyle 
a\qff =\qff \min\nolimits_{\dff\omega} U\fff(A)$\qss
and\qss
$\displaystyle
z\qff =\qff \min\nolimits_{\qff\pi} U\fff(A)$\dnsp.\oss
This implies\oss (\ref{branching-phi-image})\oss and\oss (\ref{branching-pi-image}).\oss
Since $a\qff \neq\qff z$\ffcom {\qff}it\dss follows\dss that\qss $A$\qss is\qss $\mathcal{E}$\dnsp-branching.\oss

\myitpar{Balanced almost-bases.}
An almost-basis\qss $Q\dff\in\dff\alm$\qss is said to be\qss \emph{$\omega$\dnsp-balanced 
with respect to}\qss $\mathcal{E}$\qss
if\qss $\varepsilon(Q)$\qss is also an almost-basis,\oss i.e.\qss $\varepsilon(Q)\dff\in\dff\alm$\dfcom and
\begin{equation*} 
\hspace*{1em}\varepsilon(\varphi_\omega\dff(Q))\off 
=\off \varphi_{\fff\pi\fff}(\varepsilon(Q)).
\end{equation*}
Clearly{\nsp},\qss $Q\dff\in\dff\alm$\qss 
is\qss $\pi$\nsp-balanced with respect to\qss $\mathcal{E}$\qss
if and only if\qss $\varepsilon(Q)\dff\in\dff\alm$\qss and
\begin{equation*} 
\hspace*{1em}\varepsilon(\varphi_{\fff\pi\fff}(Q))\off 
=\off \varphi_{\omega}\dff(\varepsilon(Q)).
\end{equation*}
It follows that\qss $Q\dff\in\dff\alm$\qss 
is\qss $\pi$\nsp-balanced with respect to\qss $\mathcal{E}$\qss
if and only if\qss $\varepsilon(Q)\dff\in\dff\alm$\qss and\qss $\varepsilon(Q)$\qss is\dss
$\omega$\dnsp-balanced with respect to\dss $\mathcal{E}$\dfdot
 
An almost-basis\qss $Q\dff\in\dff\alm$\qss is said to be\dff\dss \emph{$\mathcal{E}$\dnsp-balanced}\trf\qss 
if\qss $Q$\qss is both\qss $\omega$\dnsp-balanced and\qss $\pi$\dnsp-balanced with respect to\qss $\mathcal{E}$\dfdot
In view of the previous paragraph,\qss $Q\dff\in\dff\alm$\qss 
is\qss $\mathcal{E}$\dnsp-balanced if and only if\qss
$\varepsilon(Q)\dff\in\dff\alm$\qss and both\qss $Q$\qss and\qss $\varepsilon(Q)$\qss are\dss
$\omega$\dnsp-balanced with respect to\dss $\mathcal{E}$\dfdot
In particular{},\oss if\qss $Q\comm\qff\varepsilon(Q)\dff\in\dff\alm$\ffcom
then\dss $Q$\dss is\qss $\mathcal{E}$\dnsp-balanced\qss if and only if\dss 
$\varepsilon(Q)$\qss is.\qss

\mypar{Lemma.}{branching-to-balance} \emph{If\qss $A\dff\in\dff\alm$ is $\mathcal{E}$\dnsp-branching{\halfff},\oss 
then\qss $A$\qss is\qss $\mathcal{E}$\dnsp-balanced.}

\proof If\qss $A\dff\in\dff\alm$ is $\mathcal{E}$\dnsp-branching{\halfff},\oss
then\qss $\displaystyle \varepsilon\fff(\varphi_\omega(A))\off 
=\off \varepsilon\fff(A\qff +\qff a)\off
=\off A\qff +\qff z$\ffcom
and
\[
\quad
\varphi_{\fff\pi\fff}(\varepsilon\fff(A))
=\off \varphi_{\fff\pi\fff}(A)\off
=\off A\qff +\qff z\fff.
\]
It follows that\qss $A$\qss is\qss $\omega$\dnsp-balanced.\oss
By a similar argument\fff,\qss $A$\qss is\qss $\pi$\dnsp-balanced.\qss \eproof

\myitpar{Non-branching almost-bases.}
An almost-basis\dss $Q\fff\in\fff\alm$\dss is said to be\qss 
\emph{$\mathcal{E}$\dnsp-non-branching}\oss if 
\[
\varphi_\omega(Q)\off =\off \varphi_{\fff\pi\fff}(Q).
\]
In this case we will denote by\dss $\varphi\dff(Q)$\dss the coinciding images\dss
$\varphi_\omega(Q)$\dss and\dss $\varphi_{\fff\pi\fff}(Q)$\dfdot

\mypar{Lemma.}{balance-of-non-branching} \emph{Suppose that\qss $Q\comm\qff\varepsilon(Q)\in\alm$\qss
and both almost-bases\sss $Q$\sss and\sss $\varepsilon(Q)$\sss 
are\sss $\mathcal{E}$\dnsp-non-branching.\qss
Then\dss $Q$\dss is\dss $\mathcal{E}$\dnsp-balanced\qss if and only if\qss
\begin{equation}
\label{b-of-nb}
\quad
\varepsilon(\varphi\dff(Q))\qff =\qff \varphi\dff(\varepsilon(Q)),
\end{equation}
and if and only if\qss $\varepsilon(Q)$\qss is\dss $\mathcal{E}$\dnsp-balanced.}

\proof\qss The first equivalence is obvious.\qss
In order to prove the second equivalence,\qss
let us replace\qss $Q$\qss by\qss $\varepsilon(Q)$\qss in\qss (\ref{b-of-nb})\qss
and apply\qss $\varepsilon$\qss to the result.\qss
Since\qss $\varepsilon\circ\varepsilon\qff =\qff \id$\dfcom
the resulting condition is equivalent to\qss 
$\varphi\dff(\varepsilon(Q))\qff =\qff \varepsilon(\varphi\dff(Q))$\qss
and hence to\qss (\ref{b-of-nb}).\qss  \eproof

\mypar{Lemma.}{non-branching-flip} \emph{Suppose that\dss $Q$\dss is an\dss 
$\mathcal{E}$\dnsp-non-branching almost-basis and\dss
$\varphi\dff(Q)\qff =\qff Q\qff +\qff d$\nsp.\qss  
{\off}If\dff\qss $\varepsilon(Q\qff +\qff d)\qff \neq\qff Q\qff +\qff d$\dff\qss
and\dff\qss $\varepsilon(Q\qff +\qff d)\dff\in\dff\bbb$\dnsp,\qff\oss
then\qss $d\qff \neq\qff \aazz$\ffdot}

\proof\qss Suppose that\qss $d\qff =\qff a$\qss or\qss $z$\ffdot
{\qff}Then\qss $\{\dff d\comm \varepsilon(d)\trf\}\qff =\qff \{\dff\aazz\qff\}$\ffcom and
since $\varepsilon(Q\qff +\qff d) \qff \neq\qff Q\qff +\qff d$\nsp,\oss
exactly one of the elements\qss $\aazz$\qss belongs to\qss $Q\qff +\qff d$\ffdot

Therefore,\oss if\qss $d\qff =\qff a$\qss or\qss $z$\ffcom
{\off}then\qss $\aazz\dff\not\in\dff Q$\qss and hence\qss 
\[
\quad
\varepsilon(Q)\off =\off Q\fff,\quad
\varepsilon(Q\qff +\qff d)\off =\off Q\qff +\qff \varepsilon(d)\halfff.
\]
Since\qss $\varepsilon(Q\qff +\qff d)\dff\in\dff\bbb$\dfcom in this case\qss
$\varepsilon(d)\dff\in\dff U\dff(Q)$\dfcom
and hence\qss $\aazz\dff\in\dff U\dff(Q)$\dfdot

By the definition\halfff,\qss $d$\qss is the minimal element of\qss $U\dff(Q)$\qss 
with respect to both\qss $\omega$\qss and\qss $\pi$\ffdot
{\qff}It follows that if\qss $d\qff =\qff a$\sss or\sss $z$\ffcom
{\qff}then\qss $\aazz\dff\in\dff U\dff(Q)$ and one of the elements $\aazz$
is the minimal element of\qss $U\dff(Q)$\qss with respect to both\qss $\omega$\qss and\qss $\pi$\ffdot
But this contradicts to the fact that\qss 
$a\dff  <_{\dff\omega}\dff z$\qss and\qss $z\dff  <_{\dff\pi\fff}\dff a$\dfdot  \eproof
\newpage

\myitpar{Balanced bases.} A basis $B$ is said to be\dss 
\emph{$\omega$\dnsp-balanced}\dff\qss if\qss
$\varepsilon(B)\dff\in\dff\bbb$\qss and every almost-basis\qss $Q$\qss such that
\begin{equation}
\label{omega-balanced-basis}
\hspace*{1em}\mbox{either }\hspace*{0.7em}\varphi_\omega\fff(Q)\off =\off B,
\hspace*{1em}\mbox{ or }\hspace*{0.8em}
\varphi_{\fff\pi\fff}\fff(Q)\off =\off \varepsilon(B) 
\end{equation}
is $\mathcal{E}$\dnsp-balanced.\qss 
A basis is said to be\sss 
\emph{$\mathcal{E}$\dnsp-balanced}\dff\qss if\qss
it\qss is\qss \emph{both $\omega$\dnsp-balanced\qss and\qss $\pi$\dnsp-balanced}.\qss 
By interchanging the roles of\qss $\omega$\qss and\qss $\pi$\ffcom 
we see that\qss $B$\qss is\qss $\pi$\dnsp-balanced 
if and only if\qss $\varepsilon(B)\dff\in\dff\bbb$\qss and every almost-basis $Q$ such that
\begin{equation}
\label{pi-balanced-basis}
\hspace*{1em}\mbox{either }\hspace*{0.7em}\varphi_{\fff\pi\fff}\fff(Q)\off =\off B,
\hspace*{1em}\mbox{ or }\hspace*{0.8em}
\varphi_\omega\fff(Q)\off =\off \varepsilon(B)
\end{equation}
is $\mathcal{E}$\dnsp-balanced.\qss
It follows that $B$ is $\omega$\dnsp-balanced\dff\qss
if and only if $\varepsilon(B)$ is $\pi$\dnsp-balanced,\qss 
and that $B$ is $\mathcal{E}$\dnsp-balanced 
if and only if $\varepsilon(B)$ is $\mathcal{E}$\dnsp-balanced.\qss

\mypar{Lemma.}{equi-basis} \emph{If\qss $B\dff\in\dff\bbb$\dss is 
$\omega$\dnsp-balanced with respect to $\mathcal{E}$\dfcom
then\dss $\varepsilon$\dss induces a bijective map}
\begin{equation}
\label{phi-to-e-phi}
\hspace*{1em}(\varphi_\omega)^{{\minus}1}\dff(B)
\hspace*{0.8em} \longrightarrow\hspace*{0.8em} 
(\varphi_{\fff\pi\fff})^{{\minus}1}\dff(\varepsilon(B))\fff.
\end{equation}

\vspace*{-\bigskipamount}
\proof\dss If\qss $\varphi_\omega\fff (Q)\qff =\qff B$\dfcom 
then $\varepsilon(Q)\dff\in\dff\alm$\ffcom
\[
\varepsilon(\varphi_\omega\fff(Q))\off =\off \varphi_{\fff\pi\fff}\fff(\varepsilon(Q))\fff,
\]
and hence\qss $\varphi_{\fff\pi\fff}\fff(\varepsilon(Q))\qff =\qff \varepsilon(B)$\dfdot
{\qff}Similarly{\nsp},\oss if $\varphi_{\fff\pi\fff}\fff(Q)\qff =\qff \varepsilon(B)$\dfcom
{\qff}then\qss $\varepsilon(Q)\dff\in\dff\alm$\ffcom
\[
\varepsilon(\varphi_{\fff\pi\fff}\fff(Q))\off =\off \varphi_\omega\fff(\varepsilon(Q))\fff,
\]
and hence\qss $\varphi_{\omega}\fff(\varepsilon(Q))\dff =\qff \varepsilon(\varepsilon(B))\qff =\qff B$\dfdot
It follows that\qss $\varepsilon$ maps\qss $(\varphi_\omega)^{{\minus}1}\dff(B)$\qss into
$(\varphi_{\fff\pi\fff})^{{\minus}1}\dff(\varepsilon(B))$ and\dss maps\qss
$(\varphi_{\fff\pi\fff})^{{\minus}1}\dff(\varepsilon(B))$\qss into\qss
$(\varphi_\omega)^{{\minus}1}\dff(B)$\dfdot
Since $\varepsilon\circ\varepsilon=\id$\ffcom 
the two maps  induced by $\varepsilon$ are mutually inverse,\qss
and hence both of them are bijections.\qss  \eproof

\myitpar{{\dnsp}$\mathcal{E}$\dnsp-branching images.}\dss A basis $B$
is said to be an\qss \emph{$\omega$\dnsp-branching image}\oss if\oss
$B\qff =\qff \varphi_{\omega}\fff(A)$\oss
for some $\mathcal{E}$\dnsp-branching $A\dff\in\dff\alm$\dfdot
{\qff}By\qss (\ref{branching-phi-image}),\oss in this case\qss $B\qff =\qff A\qff +\qff a$\qss and\qss
\begin{equation}
\label{omega-branching-image}
\quad
\varepsilon(B)\qff =\qff A\qff +\qff z\qff =\qff \varphi_{\fff\pi\fff}\fff(A)\dff \in\dff \bbb.
\end{equation}
Trivially{\nsp},\qss $B$ is $\pi$\dnsp-branching image if and only if\qss
$B\qff =\qff \varphi_{\fff\pi\fff}\fff(A)$\qss
for some\qss $\mathcal{E}$\dnsp-branching\qss $A\dff\in\dff\alm$\dfdot
{\qff}By\qss (\ref{branching-pi-image}),\oss in this case\qss $B\qff =\qff A\qff +\qff z$\qss and\qss
\begin{equation}
\label{pi-branching-image}
\quad
\varepsilon(B)\qff =\qff A\qff +\qff a\qff =\qff \varphi_{\omega}\fff(A)\dff \in\dff \bbb.
\end{equation}
It follows that $B$ is an $\pi$\dnsp-branching image 
if and only if $\varepsilon(B)$ is an $\omega$\dnsp-branching image.\qss 
\newpage

A basis\qss $B$\qss is said to be an \emph{$\mathcal{E}$\dnsp-branching image}\qss
if it is\qss \emph{either an $\omega$\dnsp-branching image,\oss 
or\qss a $\pi$\dnsp-branching image}.\oss 
In view of the previous paragraph,\qss  
$B$\qss is an\qss 
$\mathcal{E}$\dnsp-branching image if and only if\qss 
either $B$\dfcom or $\varepsilon(B)$\qss is equal to\qss $\varphi_{\omega}\dff(A)$
for some $\mathcal{E}$\dnsp-branching\qss $A\dff\in\dff\alm$\dfdot
Therefore $B$ is an $\mathcal{E}$\dnsp-branching image if and only if $\varepsilon(B)$ is.\qss

\mypar{Lemma.}{branching-bases} \emph{If\dff\qss $B$\qss is an\qss $\mathcal{E}$\dnsp-branching image,\oss 
then\qss $\varepsilon(B)\dff \in\dff \bbb$\qss 
and\qss $\varepsilon(B)$\qss is an\qss $\mathcal{E}$\dnsp-branching image.\oss
Moreover{},\oss $\varepsilon(B)\off\neq\off B$\oss 
and\qss $B$\qss contains exactly one of the elements\qss $\aazz$\ffdot}

\proof\dss The first statement of the lemma follows from\qss 
(\ref{omega-branching-image})\qss and\qss (\ref{pi-branching-image}).\oss
If\qss $B$\qss is an $\mathcal{E}$\dnsp-branching image,\oss
then\qss $B\qff =\qff A\qff +\qff a$\qss or\qss $A\qff +\qff z$\qss for some almost-basis\qss
$A$\qss not containing\qss $\aazz$\nsp.\oss 
This implies the second statement of the lemma.\qss  \eproof

\mypar{Lemma.}{equal-fibers}  \emph{If\oss $B$\oss is not an\qss
$\mathcal{E}$\dnsp-branching image,\oss then}\oss
$\displaystyle 
(\varphi_\omega)^{{\minus}1}(B)\off =\off (\varphi_{\fff\pi\fff})^{{\minus}1}(B)$\dfdot

\vspace*{-\bigskipamount}
\proof\qss If $B$ is not an $\mathcal{E}$\dnsp-branching image,\qss
then $B$ is neither an $\omega$\dnsp-branching,\qss nor $\pi$\dnsp-branching image with respect to $\mathcal{E}$\dfdot 
Therefore,\qss if either $Q\dff\in\dff(\varphi_\omega)^{{\minus}1}(B)$\dfcom
or $Q\dff\in(\varphi_{\fff\pi\fff})^{{\minus}1}(B)$\dfcom 
then $Q$ is $\mathcal{E}$\dnsp-non-branching almost-basis,\qss and hence\qss 
$\varphi_\omega(Q)\qff =\qff \varphi_{\fff\pi\fff}(Q)$\dfdot  \eproof

\mysection{Forced\qss balance}{linking-balance}

\myitpar{The linking framework.} As before,\qss 
we assume that\qss $\omega\fff,\off\pi$\qss are two linear orders connected by an edge\qss 
$\mathcal{E}$\qss of\qss $\mathcal{L}_{{\halfff} X}$\dfcom
denote by\qss $\varepsilon$\qss the corresponding transposition\qss $\tau\dff(\mathcal{E})$\dfcom
and by\qss $\aazz$\qss be the elements interchanged by\qss $\varepsilon$\dfdot
We may assume that\qss $a\dff <_{\dff\omega}\dff z$\qss and\qss $z\dff <_{\dff\pi}\dff a$\dfdot
In the rest of this section the edge\qss $\mathcal{E}$\qss will be omitted from the notations.\qss

Let\qss $\bbb$\qss and\qss $\bbb^*$\qss be two\qss \emph{matroids}\qss on\qss $X$\ffcom
and let\qss $B\longmapsto B^*$\qss be a linking\qss  $\bbb\qff\longrightarrow\qff\bbb^*$\dnsp.\linebreak
Adjusting the notations for\qss $\bbb$\dfcom
we will denote by\qss $\alm^*$\qss the set of almost-bases of\qss $\bbb^*$\dfcom
and for\qss $Q\dff\in\dff\alm^*$\qss we will denote by\qss $U^*(Q)$\qss the set of all\qss 
$x\dff\in\dff X$\qss such that\qss $x\dff\not\in\dff Q$\qss and\qss $Q\qff +\qff x\dff\in\dff\bbb^*$\dfdot 
The same orders\qss $\omega$\dnsp,\qss $\pi$\qss and the same edge\qss $\mathcal{E}$\qss 
will be used for both pre-matroids\qss $\bbb$\qss and\qss $\bbb^*$\dfdot

\mypar{Theorem.}{linked-almost-bases}\qss \emph{Suppose that\qss $Q$\qss 
is a non-branching almost-basis of\oss $\bbb$\qss and\oss 
\begin{equation*}
\quad
\varphi\dff(Q)^*\off =\off A\qff +\qff a
\quad\mbox{ {or} }\quad
A\qff +\qff z
\end{equation*} 
for some branching almost-basis\qss $A$\qss of\qss $\bbb^*$\dfdot
Then\qss $\varepsilon\fff(Q)$\qss is also a non-branching almost-basis,\oss 
and both almost-bases\qss $Q$\qss and\qss $\varepsilon(Q)$\qss are balanced.\qss}

\proof\qss Let\qss $\displaystyle d\qff =\qff \min\nolimits_{\dff\omega} U\dff(Q)$\dfdot
{\qff}Then\qss $d\dff\in\dff U\dff(Q)$\qss and\qss 
$\displaystyle
\varphi\dff(Q)\qff =\qff Q\qff +\qff d$\ffdot
{\qff}Moreover{},\qss 
$\displaystyle (Q\qff +\qff d)^*\qff =\qff A\qff +\qff a$\qss or\qss $A\qff +\qff z$\qss
and\dss hence
\[
\quad
\varepsilon\fff((Q\qff +\qff d)^*)\dff \in\dff \bbb^*
\hspace*{0.7em}\mbox{ {and} }\hspace*{0.7em}
\varepsilon\fff((Q\qff +\qff d)^*)\qff \neq\qff (Q\qff +\qff d)^*.
\]
Now the linking condition\qss $\mathbold{L2}$\qss 
and the  injectivity of the linking map imply that\qss 
\[
\quad
\varepsilon\fff(Q\qff +\qff d)\dff \in\dff \bbb
\hspace*{0.7em}\mbox{ {and} }\hspace*{0.7em}
\varepsilon\fff(Q\qff +\qff d)\qff \neq\qff Q\qff +\qff d.
\]
Therefore we can apply Lemma\qss \ref{non-branching-flip}\qss
and conclude that\qss $d\qff \neq\qff \aazz$\ffdot

In turn{\halfff},\qss $d\qff \neq\qff \aazz$\qss implies that\qss 
$\displaystyle \varepsilon\fff(d)\qff =\qff d$ 
and hence 
$\displaystyle \varepsilon\fff(Q\qff +\qff d)\qff =\qff \varepsilon\fff(Q)\qff +\qff d$\ffdot
Since\qss 
$\displaystyle \varepsilon\fff(Q\qff +\qff d)\dff \in\dff \bbb$\dnsp,\qff\oss
it follows that\qss 
$\displaystyle \varepsilon\fff(Q)\qff +\qff d\dff \in\dff \bbb$\dnsp.\qff\oss
Therefore\qss $\varepsilon\fff(Q)\dff \in\dff \alm$\qss and\qss
$d\dff \in\dff U\dff (Q)$\dnsp.
In addition\halfff,\qss 
$\displaystyle \varepsilon\fff(Q\qff +\qff d)\qff \neq\qff Q\qff +\qff d$\qss
together with\qss $d\qff \neq\qff \aazz$\qss implies that\qss
$\varepsilon\fff(Q)\qff \neq\qff Q$\dnsp,\oss
and hence exactly one of the elements\qss $\aazz$\qss is\dss in\qss $Q$\dfdot

Suppose now that not only\dss $\displaystyle d\off =\off \min\nolimits_{\dff\omega} U\dff(Q)$\dfcom 
but also\dss 
$\displaystyle d\off =\off \min\nolimits_{\dff\omega} U\dff(\varepsilon(Q))$\dfdot
Because\qss $d\qff \neq\qff \aazz$\ffcom
{\qff}in this case the almost-basis\qss 
$\varepsilon(Q)$\qss 
is non-branching and 
\[
\quad
\varphi\dff(\varepsilon(Q))\off =\off \varepsilon(Q)\qff +\qff d.
\]
And\dss since\qss 
$\displaystyle \varepsilon(Q\qff +\qff d)\off =\off \varepsilon(Q)\dff +\dff d$\ffcom
{\qff}in this case\qss 
$\displaystyle
\varepsilon(\varphi(Q))\off 
=\off \varepsilon(Q)\qff +\qff d\off
=\off \varphi(\varepsilon(Q)).$
It follows that\qss $Q$\qss is balanced.\oss
By Lemma\qss \ref{balance-of-non-branching}\qss this implies that\qss
$\varepsilon(Q)$\qss is also balanced.\oss
It remains to prove that
$\displaystyle d\off =\off \min\nolimits_{\dff\omega} U\dff(\varepsilon(Q))$\dfdot

Let\qss $b$\qss be the element of the pair\qss $\{\dff\aazz\dff\}$\qss 
contained in\qss $Q$\dfcom and let\qss $w$\qss be the other element.\oss
Clearly{\nsp},\qss the elements\qss $b\fff,\dff\qff w$\qss are consecutive.\oss
Let\oss $C\qff =\qff Q\qff -\qff b$\ffdot
{\off}Then
\begin{equation*}
\quad
C\qff +\qff b\off =\off Q\dff\in\dff\alm\fff,
\quad 
C\qff +\qff w\off =\off \varepsilon(Q)\dff\in\dff\alm\fff,
\quad 
C\qff +\qff d\off =\off (Q\qff +\qff d)\qff -\qff b\dff\in\dff\alm\fff.
\end{equation*}

Suppose that\qss 
$\displaystyle C\qff +\qff b\qff +\qff w\qff\not\in\qff\bbb$\dnsp.\qff\oss
Then\qss $U\dff(C\qff +\qff b)\qff =\qff U\dff(C\qff +\qff w)$\qss
by Lemma\qss \ref{non-triangle}.\oss 
In other terms,\qss $U\dff(Q)\qff =\qff U\qff(\varepsilon\fff(Q))$\qss 
and hence
\begin{equation*}
\quad
d\off 
=\off \min\nolimits_{\dff\omega} U\dff(Q)\off 
=\off \min\nolimits_{\dff\omega} U\dff(\varepsilon(Q)).
\end{equation*}
It follows that in this case indeed\qss
$\displaystyle d\off =\off \min\nolimits_{\dff\omega} U\dff(\varepsilon(Q))$\dfdot  \esubproof
\vspace{0.8\bigskipamount}

\mytitle{The\qss triangular\qss case.}\dff\oss Suppose now that 
$\displaystyle C\qff +\qff b\qff +\qff w\qff \in\qff \bbb$\dnsp.\qff\qss
Since 
\begin{equation*}
\quad
C\qff +\qff b\qff +\qff d\off 
=\off Q\qff +\qff d\dff\in\dff\bbb,
\end{equation*} 
\begin{equation*}
\quad
C\qff +\qff w\qff +\qff d\off 
=\off \varepsilon(Q)\qff +\qff d\off 
=\off \varepsilon(Q\qff +\qff d)\dff \in\dff \bbb,
\end{equation*}

\vspace*{-0.7\medskipamount}
in this case $C$ together with\qss $b\fff,\dff\qff w\fff,\dff\qff d$\qss 
forms a triangle in the sense of Section \ref{matroid-triangles}.\qss

Note that in this case\qss $w\dff\in\dff U\dff(Q)$\qss because
$Q\qff +\qff w\qff =\qff C\qff +\qff b\qff +\qff w\dff \in\dff \bbb$\dfdot
This implies,\qss in particular{},\qss that $d <_{\omega} w$\nsp.\qff\oss
Since the elements $b\fff,\dff\qff w$ are consecutive,\qss 
in this case\qss $d\qff <_{\omega}\qff b$\qss also.\oss
It follows that\qss $d\qff <_{\omega}\qff a\comm\qff z$\ffdot

Let us prove that\qss 
$\displaystyle (Q\qff +\qff d)^*\off =\off A\qff +\qff b$\dnsp.\qff\oss
Recall that\qss
$\displaystyle (Q\qff +\qff d)^*\off =\off A\qff +\qff a$\oss or\oss $A\qff +\qff z$\nsp.\linebreak
In other terms,\qss
$\displaystyle (Q\qff +\qff d)^*\off =\off A\qff +\qff b$\oss or\oss $A\qff +\qff w$\nsp.\qff\oss
Suppose that\qss $\displaystyle (Q\qff +\qff d)^*\off =\off A\qff +\qff w$\dnsp.\linebreak 
Since\qss $w\dff\in\dff U\dff(Q)$\nsp,\oss
in this case one can apply Lemma\qss \ref{xy-u-sets}\qss 
to\qss $S\qff =\qff Q$\nsp,\qff\oss
$x\qff =\qff d$\nsp,\qff\oss
$y\qff =\qff w$\dnsp,\oss
and conclude that\qss $d\dff\in\dff U^*\dff(A)$\dnsp.\qff\oss
On the other hand,\qss $A$\qss is branching and hence\qss 
$\aazz$\qss are the two smallest elements of\qss $U\fff(A)$\dnsp.\oss
Since\qss $d\dff\in\dff U^*\dff(A)$\dnsp,\oss
this contradicts to\qss $d\qff <_{\omega}\qff a\comm\qff z$\nsp.\qff\oss
Therefore\qss 
$\displaystyle (Q\qff +\qff d)^*\off \neq\off A\qff +\qff w$\oss 
and\dss hence\qss
$\displaystyle (Q\qff +\qff d)^*\off =\off A\qff +\qff b$\dfdot

Let us consider an arbitrary\qss 
$\displaystyle e\dff \in\dff U\dff(\varepsilon(Q))\off 
=\off U\dff(C\qff +\qff w)$\qss
and prove that\qss
$d\qff \leq_{\omega}\qff e$\dnsp.\linebreak 
By\dss Lemma \ref{triangle}\oss either\qss
$\displaystyle e\dff\in\dff U\dff(C\qff +\qff a)\qff =\qff U\dff (Q)$\dfcom {\qff}or\oss
$\displaystyle e\dff\in\dff U\dff(C\qff +\qff d)$\dfdot
{\off}If\qss $\displaystyle e\dff\in\dff U\dff(Q)$\dfcom
then\qss $d\qff \leq_\omega\qff e$\qss by the definition of\qss $d$\dfdot
Suppose now that\qss $\displaystyle e\dff\in\dff U\dff(C\qff +\qff d)$\dnsp.\qff\oss 
Since\qss $\displaystyle d\qff <_{\omega}\qff b$\dnsp,\qff\oss
we may assume that\qss $d\qff \neq\qff e$\dnsp.\qff\oss
Note that\qss 
\[
(C\qff +\qff d)\qff +\qff b\off 
=\off (C\qff +\qff b)\qff +\qff d\off 
=\off Q\qff +\qff d
\]
and\dss hence\qss $b\dff\in\dff U\dff(C\qff +\qff d)$\dfdot 
Since,\qss as we saw above,\qss
$\displaystyle (Q\qff +\qff d)^*\off =\off A\qff +\qff b$\nsp,\qff\oss
in this case one can apply Lemma\qss \ref{xy-u-sets}\qss to\oss 
$S\qff =\qff C\qff +\qff d$\nsp,\qff\oss 
$x\qff =\qff e$\nsp,\qff\oss 
$y\qff =\qff b$\nsp,\qff\oss
and conclude that\qss $e\dff\in\dff U^*(A)$\dnsp.\qff\oss
Since\qss $\aazz$\qss are the smallest elements of $U^*(A)$\nsp,\qff\oss
it follows that $a\qff <_{\omega}\qff e$\dnsp.\qff\oss
Together with\qss $d\qff <_{\omega}\qff a$\qss this implies
that\qss $d\qff <_{\omega}\qff e$\dfdot
This completes the proof of the inequality\qss $d\qff \leq_{\omega}\qff e$\qss for all\qss 
$\displaystyle e\dff\in\dff U\dff(\varepsilon(Q))${\halfff}\dnsp.\qff\oss

It follows that in this case also\qss
$\displaystyle d\off =\off \min\nolimits_{\dff\omega} U\dff(\varepsilon(Q))$\dfdot   
The theorem follows.\qss  \esubproof   \eproof

\mypar{Corollary.}{b} \emph{Let\dff\qss $B\in\bbb$\dnsp.\oss 
If\dff\qss $B^*$\qss is an\dss $\mathcal{E}$\dnsp-branching image in $\bbb^*$\dfcom
then\qss $B$\qss is a balanced.}

\proof\qss If\dff\qss $B^*$\qss is an $\mathcal{E}$\dnsp-branching image,\oss
then\qss $\displaystyle \varepsilon(B^*)\dff \in\dff \bbb$\qss 
by Lemma\qss \ref{branching-bases}.\oss
Therefore the linking property\qss $\mathbold{L2}$\qss implies that\qss 
$\displaystyle \varepsilon(B)\dff\in\dff\bbb$\qss and\qss 
$\displaystyle \varepsilon(B^*)\off =\off \varepsilon(B)^*$\dnsp. 

It remains to prove that for every $Q\dff\in\dff\alm$
each of the two conditions\qss (\ref{omega-balanced-basis})\qss and\qss (\ref{pi-balanced-basis})\qss  
implies that $Q$ is balanced.\qss
In view of Lemma\qss \ref{branching-to-balance},\qss
we may assume that $Q$ is non-branching{\halfff}.\linebreak
In this case both conditions\qss 
(\ref{omega-balanced-basis})\qss and\qss (\ref{pi-balanced-basis})\qss mean that\qss
$\varphi\fff(Q)$\qss is equal either to\qss $B$\dfcom or\dss to\qss $\varepsilon(B)$\dfdot
It follows that\qss $\varphi\fff(Q)^*$\qss is equal either to\qss
$B^*$\dfcom or\dss 
to\qss $\varepsilon(B)^*\qff =\qff \varepsilon(B^*)$\dfdot
By Lemma\qss \ref{branching-bases}\qss $\varepsilon(B^*)$\qss 
is an\dss $\mathcal{E}$\dnsp-branching image together with\qss $B^*$\dfdot
By Theorem\qss \ref{linked-almost-bases}\qss this implies
that\qss $Q$\qss is\qss $\mathcal{E}$\dnsp-balanced.\qss  \eproof

\mypar{Lemma.}{balance-branching-images} \emph{If\qss either\qss $B$\dfcom or\qss $B^*$\qss 
is an\qss $\mathcal{E}$\dnsp-branching image,\qff\qss then 
$\varepsilon$\dss induces a bijective map}
\begin{equation*}
\quad
(\varphi_\omega)^{{\minus}1}\dff(B)\off
\longrightarrow\off 
(\varphi_{\fff\pi\fff})^{{\minus}1}\dff(\varepsilon(B)).
\end{equation*}\newpage

\vspace{-\bigskipamount}
\proof\dss By applying Corollary\qss \ref{b}\qss either to the identity linking\qss
$B\qff\longmapsto\qff B$\qss or to the linking\qss $B\qff\longmapsto\qff B^*$\qss
we see that\qss $B$\qss is balanced.\oss
It remains to apply Lemma \ref{equi-basis}.\qss  \eproof

\mysection{Coda\fff:\qff\qss the\qss order-independence}{proof}

\vspace*{\bigskipamount}
We continue to work under assumptions described 
at the beginning of Section \ref{linking-balance}.\qss

\vspace*{-0.25\bigskipamount}
\myitpar{The map\dff\qss $\bbb\qff\rightarrow\qff\bbb$\qss induced by the edge\qss $\mathcal{E}$\dnsp.}
For a basis $B\in\bbb$\dfcom
let $\varepsilon_{\dff B}\qff =\qff \varepsilon$
if either\qss $B$\dfcom or\qss $B^*$\qss is an $\mathcal{E}$\dnsp-branching image,\qss
and let\qss $\varepsilon_{\dff B}\qff =\qff \id_{\dff X}$ otherwise.\oss
Let
\[
\quad
\sigma\dff(B)\off =\off \varepsilon_{\dff B}\dff(B).
\]
By Lemma\qss \ref{branching-bases},\oss
if\qss $B\dff\in\dff\bbb$\qss is an\qss $\mathcal{E}$\dnsp-branching image,\qss
then\qss $\varepsilon(B)\dff\in\dff\bbb$\dnsp,\oss
and if\qss $B^*$\qss is an $\mathcal{E}$\dnsp-branching image,\qss
then $\varepsilon(B^*)\dff\in\dff\bbb^*$\dfdot
The linking property\qss $\mathbold{L2}$\qss implies that
in the latter case\qss $\varepsilon(B)\dff\in\dff\bbb$\qss also.\oss
It follows that\dss $\sigma$\dss is a map\qss $\bbb\qff\longrightarrow\qff\bbb$\dfdot

\vspace*{-0.25\bigskipamount}
\mypar{Lemma.}{sigma-free} \emph{$\displaystyle \sigma\dff(B)\off \neq\off B$\qss
if\qss and\qss only\qss if\qss either\qss $B$\dfcom or\qss $B^*$\qss 
is\qss an\qss $\mathcal{E}$\dnsp-branching\qss image.}\qss

\vspace*{-0.25\bigskipamount}
\proof\qss If\qss $B$\qss is an\qss $\mathcal{E}$\dnsp-branching image,\qss 
then\qss $\varepsilon(B)\qff\neq\qff B$\qss by Lemma \ref{branching-bases}.\oss
Similarly,\qss if\qss $B^*$\qss is an\qss $\mathcal{E}$\dnsp-branching image,\qss then\qss
$\varepsilon(B^*)\qff\neq\qff B^*$\dfdot
Moreover{},\qss in this case\qss
$\varepsilon(B)^*\qff =\qff \varepsilon(B^*)$\qss
by the linking property\qss $\mathbold{L2}$\qss and hence\qss
$\varepsilon(B)^*\qff\neq\qff B^*$\dfdot
By the injectivity of the linking map,\qss in this case\qss
$\varepsilon(B)\qff\neq\qff B$\qss also.\qss
Since in both cases\qss $\sigma\dff(B)\qff =\qff \varepsilon(B)$\dfcom
the\qss \emph{``if''}\qss part of the lemma follows.\qss
The\qss \emph{``only if''}\qss part follows immediately 
from the definitions.\qss  \eproof

\vspace*{-0.25\bigskipamount}
\mypar{Lemma.}{sigma-involution} \emph{If\qss
$\displaystyle \sigma\dff(B)\qff \neq\off B$\dfcom
then\qss $\varepsilon_{\dff B}\qff =\qff \varepsilon$\dnsp,\qss
$\varepsilon_{\dff \varepsilon(B)}\qff =\dff\qff \varepsilon$\dnsp,\oss
and}\qss 
\[
\quad
\sigma\dff(B)\qff =\qff \varepsilon(B)\fff,\hspace*{1.0em}
\sigma\dff(\varepsilon(B))\qff =\qff B\fff.
\]

\vspace*{-1.25\bigskipamount}
\proof\qss If\qss
$\displaystyle \sigma\off(B)\qff \neq\qff B$\dfcom then Lemma\qss \ref{sigma-free}\qss implies that
either\qss $B$\dfcom or\qss $B^*$\qss is an $\mathcal{E}$\dnsp-branching image.\oss
It follows that\qss $\varepsilon_{\dff B}\qff =\qff \varepsilon$\dnsp,\oss 
and hence\qss $\displaystyle \sigma\dff(B)\qff =\qff \varepsilon(B)$\dfdot

By Lemma\qss \ref{branching-bases},\oss
if $B$ is an $\mathcal{E}$\dnsp-branching image,\oss 
then\qss $\varepsilon(B)$\qss is also an\qss $\mathcal{E}$\dnsp-branching image.\dff\oss
Therefore,\oss in this case\qss $\varepsilon_{\dff \varepsilon(B)}\qff =\dff\qff \varepsilon$\dnsp,\oss
and hence\qss $\displaystyle \sigma\dff(\varepsilon(B))\off =\off \varepsilon(\varepsilon(B))\off =\off B$\dfdot

By Lemma\qss \ref{branching-bases}\qss applied to\qss $\bbb^*$ in the role of\qss $\bbb$\dnsp,\oss 
if\qss $B^*$\qss is and $\mathcal{E}$\dnsp-branching image,\qss 
then\qss $\varepsilon(B^*)\dff\in\dff\bbb^*$\qss and
$\varepsilon(B^*)$ is an $\mathcal{E}$\dnsp-branching image.\oss
Since\qss $\varepsilon(B^*)\dff\in\dff\bbb^*$\dnsp,\oss
the linking property\qss $\mathbold{L2}$\qss implies that\qss
$\varepsilon(B^*)\qff =\qff \varepsilon(B)^*$\dnsp.\oss  
It follows that\qss $\varepsilon(B)^*$ is an $\mathcal{E}$\dnsp-branching image.\qff\oss 
Therefore\qss $\varepsilon_{\dff \varepsilon(B)}\qff =\dff\qff \varepsilon$\qss and\qss 
$\displaystyle \sigma\dff(\varepsilon(B))\off =\off \varepsilon(\varepsilon(B))\off =\off B$\dfdot  \eproof

\mypar{Corollary.}{involution} \emph{$\sigma$\qss is an involution,\qss
i.e\qss $\displaystyle \sigma\circ\sigma\qff =\qff \id$\dnsp.\qff\oss
In particular{},\qss $\sigma$\qss is a bijection.\qss}  \eproof

\mypar{Lemma.}{sigma-invariance} \emph{For every basis\qss $B$\qss of\qss $\bbb$\dfcom
the map\qss $\varepsilon_{\dff B}$ induces bijections}
\begin{equation}
\label{phi-to-phi-sigma}
\quad
(\varphi_\omega)^{{\minus}1}\dff(B)\off\phantom{^*}
\longrightarrow\off 
(\varphi_{\fff\pi\fff})^{{\minus}1}\dff(\sigma\dff(B)),
\end{equation}
\begin{equation}
\label{phi-to-phi-sigma-star}
\quad
(\varphi_\omega^*)^{{\minus}1}\dff(B^*)\off
\longrightarrow\off 
(\varphi_{\fff\pi\fff}^*)^{{\minus}1}\dff(\sigma\dff(B^*)),
\end{equation}

\vspace{-\bigskipamount}
\proof\qss If\dss either\qss $B$\dfcom or\qss $B^*$\qss 
is\qss an\qss $\mathcal{E}$\dnsp-branching\qss image,\oss
then\qss $\displaystyle \sigma\dff(B)\off \neq\off B$\qss by Lemma\qss \ref{sigma-free},\oss
and hence\qss $\varepsilon_{\dff B}\qff =\qff \varepsilon$\qss and\qss
$\sigma\dff(B)\qff =\qff \varepsilon(B)$\qss by Lemma\qss \ref{sigma-involution}.\qff\oss
Therefore,\oss Lemma\qss \ref{balance-branching-images}\qss implies that\qss
$\varepsilon_{\dff B}$ induces a bijection\qss (\ref{phi-to-phi-sigma}).\qff\oss
If\dss neither $B$\dfcom nor $B^*$ is a branching image,\qss then
$\varepsilon_{\dff B}\qff =\qff \id_{\dff X}$ 
and $\sigma\dff(B)\qff =\qff B$\dnsp.\qff\oss
Therefore,\oss in this case Lemma\qss \ref{equal-fibers}\qss implies that
$\varepsilon_{\dff B}$ induces a bijection\qss (\ref{phi-to-phi-sigma}).\qss

Note that replacing the linking\qss $\bbb\qff\longrightarrow\qff\bbb^*$\qss
by the inverse linking\qss $\bbb^*\qff\longrightarrow\qff\bbb$\qss
does not affect neither\qss $\varepsilon_B$\nsp,\dff\oss nor\qss $\sigma$\dnsp.\qff\oss
Therefore,\qss by applying the already proved part of the lemma 
to the inverse linking\qss $\bbb^*\qff\longrightarrow\qff\bbb$\dnsp,\oss
we see that\qss $\varepsilon_{\dff B}$ 
induces a bijection\qss (\ref{phi-to-phi-sigma-star}).\qss  \eproof

\mypar{Corollary{\nsp}.}{card-omega} \emph{For every basis\qss $B$\qss of\qff\qss $\bbb$\qss
the following equalities hold.}
\begin{equation*}
\quad
\card\dff (\varphi_\omega)^{{\minus}1}\dff(B)\phantom{^*}
\off =\off 
\card\dff (\varphi_{\fff\pi\fff})^{{\minus}1}\dff(\sigma\dff(B)).
\end{equation*}
\begin{equation*}
\quad
\card\dff (\varphi_\omega^*)^{{\minus}1}\dff(B^*)
\off =\off 
\card\dff (\varphi_{\fff\pi\fff}^*)^{{\minus}1}\dff(\sigma\dff(B)^*).\qss  \mbox{ \eproof }
\end{equation*}

Let us restate the last corollary in terms of the Tutte activities,\qss
and then in terms of the maps\qss $\mathfrak{W}_{\dff\bullet}\dff({\dff\bbb\qff\to\qff\bbb^*})$\qss 
(see Section \ref{orders-pre-matroids}).\qss

\mypar{Corollary{\nsp}.}{activity-omega} $\displaystyle
i_{\dff\omega}\dff(B)\off =\off 
i_{\dff\pi\fff}\dff(\sigma\dff(B))$\qss
\emph{and}\qff\qss
$\displaystyle 
i_{\dff\omega}\dff(B)\off =\off 
i_{\dff\pi\fff}\dff(\sigma\dff(B)^*)$\qss
\emph{{\qff}for every\qss $B\dff\in\dff\bbb$\dfdot}  \eproof

\mypar{Corollary{\nsp}.}{w-maps-sigma}
$\displaystyle 
\mathfrak{W}_{\dff\omega}\dff({\dff\bbb\qff\to\qff\bbb^*}) 
\off =\off
\mathfrak{W}_{\dff\pi\fff}\dff({\dff\bbb\qff\to\qff\bbb^*})\circ\sigma$\dfdot  \eproof

\myitpar{Proof\qss of\qss Theorem\qss \ref{w-independence}.} 
Suppose now that $\omega$ and $\pi$ are two arbitrary linear orders on $X$\dfdot
In view of Lemma \ref{multi-image-composition},\qss
the last corollary implies that
\begin{equation}
\label{omega-pi}
\quad
\mathbb{W}_{\dff\omega}\dff({\dff\bbb\qff\to\qff\bbb^*}) 
\off =\off
\mathbb{W}_{\dff\pi\fff}\dff({\dff\bbb\qff\to\qff\bbb^*})
\end{equation}
if the orders $\omega$ and $\pi$ are  connected by an edge of\qss $\mathcal{L}_{{\halfff} X}$\dfdot
In view of Lemma\qss \ref{connecting-pairs},\qss
the equality\qss (\ref{omega-pi})\qss for linear orders $\omega\comm\pi$ connected by an edge
implies the equality\qss (\ref{omega-pi})\qss for arbitrary two linear orders $\omega\comm\pi$\dfdot
It follows that the Whitney multi-set
$\displaystyle \mathbb{W}_{\dff\omega}\dff({\dff\bbb\qff\to\qff\bbb^*})$\qss
of\qss $\bbb\qff\to\qff\bbb^*$\qss with respect to\qss $\omega$\qss 
does not depend on the choice of\qss $\omega$\dfdot  \eproof

\myappend{The\qss symmetric\qss exchange\qss property}{symmetric exchange}

\vspace*{\bigskipamount}
In this appendix\qss $\bbb$\qss is assumed to be a\qss \emph{matroid}\qss on\qss $X$\dfdot

\vspace*{-0.25\bigskipamount}
\myappar{Lemma.}{no-proper} \emph{No basis of\fff\qss $\bbb$\qss is properly contained in another basis.\oss
In particular{},\oss no almost-basis is a basis and no over-basis is a basis.}

\vspace*{-0.25\bigskipamount}
\proof\qss If\qss $B_1\fff,\off B_2\qff \in\dff \bbb$\dfcom $B_1\qff \subset\qff B_2$\ffcom
and\qss $B_1\qff \neq\qff B_2$\nsp,\oss
then\qss $B_2 \smallsetminus B_1\qff \neq\qff \varnothing$\dfdot
By applying the exchange property to any\qss $x\dff \in\dff B_2 \smallsetminus B_1$\ffcom
we see,\qss in particular{},\qss that\qss $B_1 \smallsetminus B_2$\qss is non-empty\nsp.\oss
The contradiction with\qss $B_1\qff \subset\qff B_2$\qss proves the first statement of the lemma,\qss
which immediately implies the second one.\qss  \eproof

\vspace*{-0.25\bigskipamount}
\myappar{Lemma.}{max-over} \emph{Suppose that\qss $J\qff\subset\qff S\subset X$\dfdot
If\qss $S$\qss contains a basis and\qss $J$\qss is contained in a basis,\qss
then there exist a basis\qss $B$\qss such that\qss $J\qff\subset\qff B\qff\subset\qff S$\dfdot}

\vspace*{-0.25\bigskipamount}
\proof\qss Consider all bases containing $J$ and choose among them a basis $B$ such that 
$B\cap S$ is maximal with respect to the inclusion.\dss
Note that $J\qff\subset\qff B\cap S$ because $J\qff\subset\qff B$ and $J\qff\subset\qff S$\dfdot
It is sufficient to prove that $B\qff\subset\qff S$\dfdot

Suppose that $B\dff\not\subset\dff S$\dfdot
Then $B\dff\smallsetminus\dff S\qff\neq\qff\varnothing$\dfdot
Let $x\dff\in\dff B\dff\smallsetminus\dff S$\dfdot
Consider some basis $\beta\qff\subset\qff S$\dfdot
Clearly{}, $x\dff\in\dff B\dff\smallsetminus\dff\beta$\dfdot
By the exchange property{}, $B\qff -\qff x\qff +\qff y\dff\in\dff\bbb$ 
for some $y\dff\in\dff\beta\dff\smallsetminus\dff B$\dfdot
Since $x\dff\not\in\dff S$ and,\dss obviously{},\qss 
$y\dff\in\dff S\dff\smallsetminus\dff B$\dfcom we see that
\begin{equation*}
\quad
(B\qff -\qff x\qff +\qff y)\cap S\off =\off (B\cap S)\qff +\qff y\fff.
\end{equation*}
and hence $(B\qff -\qff x\qff +y)\cap S$ properly contains $B\cap S$\dfcom
and at the same time $J\qff\subset\qff B\cap S\qff\subset\qff (B\qff -\qff x\qff +\qff y)\cap S$\dfdot
But this contradicts to the choice of $B$\dfdot
It follows that $B\qff\subset\qff S$\dfdot  \eproof

\myitpar{Independent sets and circuits.} A subset $Y\qff\subset\qff X$ is called\qss \emph{independent}\qss 
if\qss $Y\qff\subset\qff B$\qss for some\qss $B\dff\in\dff\bbb$\dfcom
and\qss \emph{dependent}\qss if\qss $Y$\qss is not contained in any\qss $B\in\bbb$\dfdot
Clearly\nsp,\qss any subset of an independent set is independent.\dss
A subset\qss $C\subset X$\qss is called a\qss \emph{circuit}\dff\qss if\qss $C$\qss is dependent,\qss
but\qss $C\qff -\qff x$\qss is independent for every\qss $x\dff\in\dff C$\dfdot

\vspace*{-0.25\bigskipamount}
\myappar{Lemma.}{curcuit-of-over-basis} \emph{For every over-basis\qss $Q$\qss 
the set\qss $C\dff(Q)$\qss is a circuit.}

\vspace*{-0.25\bigskipamount}
\proof\qss By the definition, $x\dff\in\dff C\dff(Q)$ if and only if $Q\qff -\qff x\dff\in\dff\bbb$\dfdot
Therefore,\qss if $x\dff\in\dff C\dff(Q)$\dfcom then $C(Q)\qff -\qff x\qff\subset\qff Q\qff -\qff x\dff\in\dff\bbb$\dfdot
It follows that $C(Q)-x$ is independent for every $x\dff\in\dff C\dff(Q)$\dfdot
It remains to prove that $C(Q)$ is dependent.\qss

Suppose that $C\dff(Q)$ is independent.\qss 
Then $C\dff(Q)$ is contained in some basis.\qss
Since $Q$ is an over-basis,\dss $Q$ contains a basis.\qss
By applying Lemma \ref{max-over} to $J\qff =\qff C\dff(Q)$ and $S\qff =\qff Q$\dfcom
we see that $C\dff(Q)\qff \subset\qff  B\qff \subset\qff  Q$ for some $B\dff\in\dff\bbb$\dfdot
By Lemma\qss \ref{no-proper}\qss $Q$\qss is not a basis,\qss
and hence\qss $B\neq Q$\dnsp.\oss 
It follows that\qss $Q\dff\smallsetminus\dff B\dff\neq\dff\varnothing$\dfdot

Let $q\dff\in\dff C\dff(Q)$\dfdot
Then $q\dff\in\dff B$ and $Q\qff -\qff q\dff\in\dff\bbb$\dfdot
Consider some $x\dff\in\dff Q\dff\smallsetminus\dff B\qff =\qff (Q\qff -\qff q)\dff\smallsetminus\dff B$\dfdot
By the exchange property{},\qss 
\[
\quad
(Q\qff -\qff q)\qff -\qff x\qff +\qff y\qff\in\qff\bbb
\]
for some $y\in B\dff\smallsetminus\dff (Q\qff -\qff q)$\dfdot
But $B\dff\smallsetminus\dff (Q\qff -\qff q)\qff =\qff \{\dff q\dff\}$ 
because\qss $B\qff \subset\qff Q$\qss and\qss $q\dff\in\dff B$\dfdot
It follows that $y\qff =\qff q$ and hence 
\[
\quad
(Q\qff -\qff q)\qff -\qff x\qff +\qff y\off 
=\off Q\qff -\qff q\qff -\qff x\qff +\qff q\off 
=\off Q\qff -\qff x\fff.
\]
It follows that $Q\qff -\qff x\dff\in\dff\bbb$ and hence $x\dff\in\dff C\dff(Q)$\dfdot
But $x\dff\not\in\dff B$ and hence $x\dff\not\in\dff C\dff(Q)\subset B$\dfdot
Therefore,\dss the assumption that $C\dff(Q)$ is independent leads to a contradiction.\qss
It follows that $C\dff(Q)$ is a dependent set.\qss  \eproof

\myappar{Lemma.}{no-basis} \emph{Let\qss $A$\qss be an almost-basis,\dss
and let\qss $Y\qff \subset\qff  X$\qss be a subset disjoint from\qss $A$\dfdot
If\qss $A\qff +\qff y\dff\not\in\dff\bbb$\qss for all\qss $y\dff\in\dff Y$\dfcom
then\qss $A\cup Y$\qss does not contain any basis.\dss}

\proof\qss Since $A$ is an almost-basis,\dss $A\qff +\qff a\in\bbb$ for some $a\dff\not\in\dff A$\dfdot
Then,\qss in particular, $a\dff\not\in\dff Y$\dfcom
and hence $a\dff\not\in\dff A\cup Y$\dfdot

Suppose that $B\dff\in\dff\bbb$ and $B\qff \subset\qff  A\cup Y$\dfdot
Then $a\dff\not\in\dff B$ and hence
\begin{equation*}
\quad
B\dff\smallsetminus\dff (A\qff +\qff a)\off 
=\off B\dff\smallsetminus\dff A\fff,\hspace*{1.4em} 
(A\qff +\qff a)\smallsetminus B\off 
=\off (A\dff\smallsetminus\dff B)\qff +\qff a\fff.
\end{equation*}
In particular,\dss $a\dff\in\dff (A\qff +\qff a)\dff\smallsetminus\dff B$\dfdot
By the exchange property of matroids,\dss
there exists an element $y\dff\in\dff B\dff\smallsetminus\dff (A\qff +\qff a)\qff =\qff B\dff\smallsetminus\dff A$
such that $(A\qff +\qff a)\qff -\qff a\qff +\qff y\qff =\qff A\qff +\qff y$ is a basis.\qss
But $B\dff\smallsetminus\dff A\qff \subset\qff  Y$ and hence $y\dff\in\dff Y$\dfdot
This contradicts to the assumption that $A\qff +\qff y$ is not a basis for all $y\dff\in\dff Y$\dfdot
It follows that $A\cup Y$ cannot contain a basis.\dss  \eproof

\myappar{Lemma.}{circuit-no-basis} \emph{Let\qss $A$\qss be an almost-basis,\qss
and let\qss $Y\qff \subset\qff  X$\qss be a subset disjoint from\qss $A$\dfdot
Suppose that\qss $A\qff +\qff y\dff\not\in\dff\bbb$\qss for all\qss $y\dff\in\dff Y$\dfdot
If\qss $C$\qss is a circuit such that\qss $C\qff -\qff c\qff \subset\qff  A\cup Y$\qss 
for some\qss $c\dff\in\dff C$\dfcom $c\dff\not\in\dff A\cup Y$\dfcom
then\qss $A\cup (Y\qff +\qff c)$\qss does not contain any basis.\dss}

\proof\qss Suppose that $A\cup (Y\qff +\qff c)$ contains a basis.\qss
Since $C$ is a circuit,\dss the set $C\qff -\qff c$ is independent,\qss 
i.e. is contained in a basis.\qss
By applying Lemma \ref{max-over} to $J\qff =\qff C\qff -\qff c$ and $S\qff =\qff A\cup (Y\qff +\qff c)$\dfcom
we see that $C\qff -c\qff \subset\qff  B\qff \subset\qff  A\cup (Y\qff +\qff c)$ for some $B\dff\in\dff\bbb$\dfdot

If $c\dff\in\dff B$\dfcom then 
$C\qff =\qff (C\qff -\qff c)\qff +\qff c\qff \subset\qff  B$ 
and hence $C$ is independent.\qss
Since $C$ is a circuit,\dss this is impossible.\dss
It follows that $c\not\in B$ and hence $B\qff \subset\qff  A\cup Y$\dfdot
But the last inclusion contradicts to Lemma \ref{no-basis}.\qss
Therefore,\dss $A\cup (Y\qff +\qff c)$ cannot contain a basis.\qss
This proves the lemma.\qss  \eproof

\myappar{Lemma.}{u-c-intersection} \emph{Suppose that\qss $A$\qss 
is an almost-basis and\qss $Q$\qss is an over-basis.\oss 
If the intersection\qss $U\dff(A)\cap C\dff(Q)$\qss is not empty,\qss 
then it contains at least\qss $2$\qss elements.\dss}

\proof\qss Suppose that $U\dff(A)\cap C\dff(Q)$ consists of only one element $c$\dfdot
Then $A\qff +\qff y\dff\not\in\dff\bbb$ for all $y\dff\in\dff C\dff(Q)\qff -\qff c$\dfcom $y\dff\not\in\dff A$\dfdot 
Let $C\qff =\qff C\dff(Q)$ and $Y\qff =\qff (C\dff\smallsetminus\dff A)\qff -\qff c$\dfdot
Then $A$\dfcom $Y$\dfcom $C$\dfcom and $c$ satisfy the assumptions of Lemma \ref{circuit-no-basis}.\qss
This lemma implies that $A\cup(Y\qff +\qff c)$ does not contain any basis.\qss
But $Y\qff +\qff c\qff =\qff C\dff\smallsetminus\dff A$ and hence $A\cup(Y\qff +\qff c)\qff =\qff A\cup C$\dfdot
It follows that $c\dff\in\dff A\cup(Y\qff +\qff c)$ and hence $A\cup(Y\qff +\qff c)$ contains $A\qff +\qff c$\dfdot
But $c\dff\in\dff U\dff(A)$ and hence $A\qff +\qff c$ is a basis.\qss
The contradiction shows that $U\dff(A)\cap C\dff(Q)$ cannot consist of only one element.\qss
The lemma follows.\qss  \eproof

\myitpar{Proof of Theorem \ref{symm-exchange}.} Suppose that $B_1\dff,\off B_2\in\bbb$ 
and $x\dff\in\dff B_1\dff\smallsetminus\dff B_2$\dfdot
Let 
\begin{equation*} 
\hspace*{1.0em}
E_1\off =\off \{\dff y\dff \in\dff  B_2\dff \smallsetminus\dff  B_1\dff \mid\dff  
B_2\qff -\qff y\qff +\qff x\dff\in\dff\bbb\dff \}
\end{equation*}
\begin{equation*}
\hspace*{1.0em}
E_2\off =\off \{\dff y\dff \in\dff  B_2\dff \smallsetminus\dff  B_1\dff \mid\dff  
B_1\qff -\qff x\qff +\qff y\dff\in\dff\bbb\dff \}
\end{equation*}

\vspace*{-\medskipamount}
It is sufficient to prove that $E_1\cap E_2\qff\neq\qff\varnothing$\dfdot

Note that $B_2+x$ is an over-basis,\dss and the condition 
$B_2\qff -\qff y\qff +\qff x\dff\in\dff\bbb$
is equivalent to 
$(B_2\qff +\qff x)\qff -\qff y\dff\in\dff\bbb$\dfdot
Note also that 
$(B_2\qff +\qff x)\dff\smallsetminus\dff B_1\qff =\qff B_2\dff\smallsetminus\dff B_1$ 
because $x\dff\in\dff B_1$\dfdot
It follows that $E_1\qff =\qff C\dff(B_2\qff +\qff x)\dff\smallsetminus\dff B_1$\dfdot
Similarly{}, $B_1\qff -\qff x$ is an almost-basis,\dss 
and the condition $B_1\qff -\qff x\qff +\qff y\dff\in\dff\bbb$
is equivalent to $(B_1\qff -\qff x)\qff +\qff y\dff\in\dff\bbb$\dfdot
Note also that $(B_1\qff -\qff x)\cap B_2\qff =\qff B_1\cap B_2$ 
because $x\dff\not\in\dff B_2$\dfdot
It follows that $E_2\qff =\qff U\dff(B_1\qff -\qff x)\cap B_2$\dfdot

Note that $x\dff\in\dff C\dff(B_2\qff +\qff x)$ and $x\dff\in\dff U\dff(B_1\qff -\qff x)$\dfcom
and hence the intersection $U\dff(B_1\qff -\qff x)\cap C\dff(B_2\qff +\qff x)$ is not empty.\qss
By Lemma\qss \ref{u-c-intersection}\qss it consists of at least $2$ elements.\qss
Let $y$ be some element of this intersection different from $x$\dfdot

Then $y\dff\not\in\dff B_1$ because $y\qff \neq\qff  x$ and 
$U\dff(B_1\qff -\qff x)$ is disjoint from $B_1\qff -\qff x$\dfdot
Similarly, $y\dff\in\dff B_2$ because $C\dff(B_2\qff +\qff x)\qff \subset\qff  B_2\qff +\qff x$ and $y\qff \neq\qff  x$\dfdot
It follows that $y\dff\in\dff E_1\qff  =\qff  C\dff(B_2\qff +\qff x)\dff\smallsetminus\dff B_1$ 
and $y\dff\in\dff E_2\qff  =\qff  U\dff(B_1\qff -\qff x)\cap B_2$\dfdot
Therefore $y\dff\in\dff E_1\cap E_2$ and hence $E_1\cap E_2\qff \neq\qff \varnothing$\dfdot
The theorem follows.\dss  \eproof

\myappar{Remark.}{fundamental-curcuit} \emph{For every over-basis\qss $Q$\qss 
the set\qss $C\dff(Q)$\qss is the only circuit contained in\qss $Q$.\dfdot}

\proof Let $C\qff \subset\qff  Q$ be an arbitrary circuit in $Q$\dfdot
If $x\dff\in\dff C\dff(Q)$\dfcom but $x\dff\not\in\dff C$\dfcom
then $C\subset Q\qff -\qff x$ and $Q\qff -\qff x$ is a basis.\qss
In this case $C$ is independent,\qss 
in contradiction with being a circuit.\qss
It follows that $C\dff(Q)\qff \subset\qff  C$\dfdot

If $C\dff(Q)\qff \neq\qff  C$\dfcom then $C\dff(Q)\qff \subset\qff  C$ imlies that
$C\dff(Q)\qff \subset\qff  C\qff -\qff x$ for some $x\dff\in\dff C$\dfdot 
Since $C$ is a circuit,\qss $C\qff -\qff x$ is independent,\qss
and hence $C\dff(Q)$ is also independent.\qss 
But $C\dff(Q)$ is a circuit by Lemma \ref{curcuit-of-over-basis},\qss
and hence is a dependent set.\qss
It follows that $C\qff =\qff C\dff(Q)$\dfdot  \eproof

\myappend{Duality}{duality}

\myitpar{The dual of a pre-matroid.} As usual,\qss let\qss $\bbb$\qss be a pre-matroid on $X$\dfcom
and let\qss $\cbbb$\qss be its dual pre-matroid\qss (see Section \ref{basic}).\qss
Let\qss $\alm\csup$\qss and\qss $\ove\csup$\qss be,\qss respectively{\nsp},\qss 
the sets of all almost-bases and all over-bases of\qss $\bbb\csup$\dfdot 
For\qss $D\dff\in\dff\alm\csup$\qss we will denote by\qss $U\csup\dff(D)$\qss 
the set of all\qss $x\dff\in\dff X$\qss 
such that\qss $x\dff\not\in\dff D$\qss and\qss $D\qff +\qff x\dff\in\dff\bbb\csup$\dfdot
Similarly{\nsp},\qss for\qss $Q\dff\in\dff\ove\csup$\qss 
we will denote by\qss $C\dff\csup(Q)$\qss the set of all\qss $x\dff\in\dff Q$
such that $Q\qff -\qff z\dff\in\dff\bbb\csup$\dfdot

\myitpar{Proof of Theorem \ref{dual-matroid-1}.} Let $\bbb$ be a matroid.\dss
By Theorem \ref{symm-exchange} $\bbb$ has the symmetric exchange property{}.\qss
By taking complements the symmetric exchange property of $B$
translates into the following property of $\cbbb$\dfdot
\vspace*{-\bigskipamount}\vspace*{\medskipamount} 
\begin{quote}
\emph{If\dff\qss $B_1\fff,\off B_2$\qss 
are bases of\dff\qss $\cbbb$\qss and\qss $x\dff\in\dff B_2\smallsetminus B_1$\ffcom 
then there exists $y\dff\in\dff B_1\smallsetminus B_2$ \\ 
such that both\qss $B_2\dff\sdiff\dff\{\dff x\comm y\dff\}$\qss 
and\qss $B_1\dff\sdiff\dff\{\dff x\comm y\dff\}$\dss are bases of\qss $\cbbb$\nsp.}
\end{quote}
\vspace*{-\bigskipamount}\vspace*{\medskipamount}
This property turns into the symmetric exchange property of $\cbbb$
after interchanging the roles of $B_1$ and $B_2$\dfdot
It follows that the dual $\cbbb$ of $\bbb$ also satisfies the symmetric exchange property{}.\dss
Since the symmetric exchange property trivially implies the exchange property{},\dss
the dual $\cbbb$ has the exchange property{}.\dss  
Therefore $\cbbb$ is a matroid.\dss  \eproof

\myitpar{Remark.} The symmetric exchange property cannot be replaced in the proof of Theorem \ref{dual-matroid-1} by the exchange property{}.\dss
Indeed,\dss taking complements turns the exchange property of\qss $\bbb$\qss into the following property of\qss $\cbbb$\dfdot
\vspace*{-\bigskipamount}\vspace*{\medskipamount} 
\begin{quote}
\emph{If\dff\qss $B_1\fff,\off B_2$\qss 
are bases of\qss $\cbbb$\qss and\qss $x\dff\in\dff B_2\smallsetminus B_1$\ffcom then\qss \\ 
$B_1\qff +\qff x\qff -\qff y$\dss is a basis of\qss $\cbbb$\qss for some\qss $y\dff\in\dff B_1\smallsetminus B_2$\nsp.}
\end{quote}
\vspace*{-\bigskipamount}\vspace*{\medskipamount}
This property is different from the exchange property{},
and does not turns into the exchange property after interchanging the roles of $B_1$ and $B_2$\dfdot

\myappar{Lemma}{almost-and-over-bases} \emph{Let\qss $Q\subset X${\halfff}\dfdot
Then\qss $Q\dff\in\dff\ove$\qss if and only if\qss $Q\csup\dff\in\dff\alm\csup$\dfdot
If\qss $Q\in\ove$\dfcom then}
\begin{equation}
\label{c-u-duality}
\hspace*{1em}
C\dff(Q)\csup\off =\off U\csup\dff(Q\csup).\qquad\mbox{ \eproof }
\end{equation}

\vspace*{-\bigskipamount}
\myitpar{Maps defined by a linear order.} Suppose that a linear order $\omega$ on $X$ is fixed.\qss
Let 
\[
\hspace*{1em}
\varphi_{\omega}^{\ccomp}\dff\colon\dff\alm\csup\off\longrightarrow\off\bbb\csup\nsp,
\qquad\psi_{\omega}^{\ccomp}\dff\colon\dff\ove\csup\off\longrightarrow\off\bbb\csup
\]
be the maps defined as the maps $\varphi_{\omega},\qff\off \psi_{\omega}$
defined in Section\qss \ref{orders-pre-matroids},\qff\qss 
but with $\bbb\csup$ playing the role of $\bbb$\dfdot
For\qss $B\in\bbb\csup$\qss let\oss $i_{\dff\omega}^{\ccomp}\dff(B)$\oss
be the number of elements in\qss $(\varphi_{\omega})^{{\minus}1}\fff(B)$\dnsp.\oss 
Lemma\qss \ref{almost-and-over-bases}\qss immediately implies the following lemma.

\myappar{Lemma.}{phi-psi-dual} \emph{If\qss $Q\dff\in\dff\ove$\dfcom then}
$\displaystyle \psi_{\omega}\dff(Q)\csup\off =\off \varphi_{\omega}^{\ccomp}\trf(Q\csup\fff)$\dfdot\quad   \eproof

\myappar{Lemma.}{alpha-omega-dual} \emph{If\qff\qss $B\in\bbb$\dfcom then}\qss
$\displaystyle e_{\dff\omega}\dff(B)\off =\off 
i_{\dff\omega}^{\ccomp}\dff(B\csup{\halfff})$\dfdot

\proof\qss Let\qss $B\in\bbb$\dfdot
Lemma\qss \ref{phi-psi-dual}\qss implies that\qss $B\qff = \psi\dff(Q)$ 
if and only if\qss $B\fff\csup\qff =\qff \varphi_{\omega}^{\ccomp}\dff(Q\fff\csup)$\dnsp.\qff\oss
Therefore,\qss the map\qss $Q\qff\longmapsto\qff Q\csup$\qss induces a bijection
\[
\hspace*{1em}
\psi^{\dff -1}\dff(B)\off \longrightarrow\off (\varphi_{\omega}^{\ccomp})^{\dff -1}\dff(B\csup),
\]
and hence these two sets have the same number of elements.\oss
The lemma follows.\qss  \eproof

\myappend{Permutations\qss and\qss transpositions}{perm-trans}

\myappar{Lemma.}{permutations-generators} \emph{Suppose that\qss $n\dff\in\dff\nnn$\qss and\qss $n\qff \geq\qff 2$\dnsp.\oss
Let\qss $X_{\dff n}\qff =\qff \{\dff 1\comm\qff 2\comm\qff \ldots\dff,\dff\qff n\dff\}${\halfff}\dfcom
and for each\oss $i\qff =\qff 1\comm\qff 2\comm\qff \ldots\dff,\dff\qff n-1$\oss 
let\qss $\varepsilon_i$\qss be the transposition of elements\qss $i\comm\qff i+1$\dnsp.\qff\oss
Then every permutation of\qss $X_{\dff n}$\qss is equal to a composition of 
several transpositions of the form\qss $\varepsilon_i$\ffdot}

\proof\qss We use induction by $n$\dfdot
There are only two permutations of $X_{\dff 2}$\ffcom namely{},\oss 
$\id_X$\qss and\qss $\varepsilon_1$\dfdot
Let\qss $m\dff \in\dff \nnn$\qss and\qss $m\qff >\qff 2$\nsp.\oss
Suppose that the lemma is true for all $n<m$\dnsp,\oss
and consider a permutation\qss $\sigma$\qss of the set\qss 
$X_{\dff m}\qff =\qff \{\dff 1\comm\qff 2\comm\qff \ldots\dff,\dff\qff m\dff\}$\dfdot

If\qss $\sigma(m)\qff =\qff m$\dfcom then\qss $\sigma$\qss 
induces a permutation of\qss $X_{\dff m-1}$\dfdot
By the inductive assumption the induced permutation 
is equal to a composition of several transpositions 
of\qss $X_{\dff m-1}$\qss of the required form.\dss
Then\qss $\sigma$\qss is equal to the composition of transpositions of the same elements,\dss
but considered as elements of\qss $X_m$\qss (note that all of them leave\qss $m$\qss fixed).\dss

It remains to consider the case when $\sigma\dff(m)<m$\dfdot
Let $a\qff =\qff \sigma\dff(m)$\dfdot
Then 
\[
\quad
\varepsilon_{m\dff -\dff 1}\circ\ldots\circ\varepsilon_{a\dff +\dff 1}\circ\varepsilon_a\circ\sigma\dff(m)\off =\off m\dff.
\]
By the previous paragraph the permutation $\tau\qff =\qff \varepsilon_{m-1}\circ\ldots\circ\varepsilon_a\circ\sigma$
is equal to a composition of transpositions $\varepsilon_i$\ffdot
Since every transposition is equal to its own inverse,\dss
\[
\quad
\sigma\off =\off \varepsilon_a\circ\varepsilon_{a+1}\circ\ldots\circ\varepsilon_{m-1}\circ\tau\dff,
\]
and hence $\sigma$ is also equal to a composition of transpositions $\varepsilon_i$\ffdot  \eproof

\prooftitle{Proof of Lemma \ref{connecting-pairs}}  We may assume that\qss 
$X\qff =\qff X_{\dff n}\qff =\qff \{\fff 1\dff,\dff 2\dff,\dff \ldots\dff,\dff n\fff\}$\qss 
for some\qss $n\dff\in\dff\nnn$\dnsp,\oss
and that the order\qss $\omega$\qss is induced by the standard order on $\nnn$\dnsp.\oss
Any other linear order on\dss $X$\dss has the form\dss $\sigma\cdot\omega$\dss 
for some permutation\dss $\sigma$\dss of\dss $X$\dfdot

Let\qss $\varepsilon_{\dff i}$\qss be the transpositions defined in Lemma\qss \ref{permutations-generators}.\oss
If $\tau$ is a permutation of $X$\ffcom
then $\tau_{\dff i}=\tau\circ\varepsilon_{\dff i}\circ\tau^{{\minus}1}$ 
is the transposition of elements $\tau\dff(i)$ and $\tau\dff(i+1)$\dfcom
which are consecutive with respect to the order\qss $\tau\cdot\omega$\dnsp.\oss
At the same time\qss $\tau_{\dff i}\circ\tau\qff =\qff \tau\circ\varepsilon_{\dff i}$\qss and hence 
\begin{equation}
\label{orders-transp}
\quad
\tau_{\dff i}\cdot(\dff\tau\fff\cdot\omega)\off 
=\off (\tau_{\dff i}\fff\circ\fff\tau\fff)\cdot\omega
=\off (\fff\tau\dff\circ\dff\varepsilon_{\dff i})\cdot\omega.
\end{equation}
By Lemma \ref{permutations-generators} every permutation $\sigma$ of $X$ is equal to a composition of several transpositions of the form $\varepsilon_{\dff i}$\ffdot
It follows that one can get from $\id_X$ to $\sigma$ by a sequence of right compositions with $\varepsilon_i$\ffdot
In view of\qss (\ref{orders-transp}),\qff\qss this implies that one can connect the standard order $\omega$ 
with the order $\sigma\cdot\omega$ by a sequence of orders of the form 
\[
\quad
\omega\off =\off \tau_{\dff 0}\cdot\omega\fff,\off\off\off  
\tau_{\dff 1}\cdot\omega\fff,\off\off\off  
\tau_{\dff 2}\cdot\omega\fff,\off\off\off   
\ldots\fff,\off\off\off  
\tau_{\dff m}\cdot\omega\off   
=\off \sigma\cdot\omega
\]
(where,\qss of course,\qss $\tau_{\dff 0}\qff =\qff \id$\qss and\qss $\tau_{\dff m}\qff =\qff \sigma$\nsp)\qss
such that $\tau_{\dff i\dff +\dff 1}\cdot\omega$ is the image of $\tau_{\dff i}\cdot\omega$
under the transposition of two elements consecutive with respect to $\tau_{\dff i}\cdot\omega$
for all $i=0\comm\qff 1\comm\qff 2\comm\qff\ldots\comm\qff m-1$\dfdot
The lemma follows.\qss  \eproof

\myappend{A\qss direct\qss proof\qss of\qss Lemma\qss \ref{xy-u-sets}}{linkings-lemma-proof}

\vspace*{\bigskipamount}
Let\qss $\bbb$\dfcom $\bbb^*$\qss be pre-matroids on\qss $X$\ffcom
and let\dff\qss $L\dff\colon \bullet\qff\longmapsto\qff\bullet^*$\dff\qss 
be a linking\qss $\bbb\qff\longmapsto\qff\bbb^*$\dfdot

\myappar{Lemma.}{xy-flip}\qss \emph{Suppose that\dss $S$\dss is an almost-basis of\qss $\bbb$\dss and\qss
$x\comm y\dff\in\dff U\dff(S)$\ffcom {\qff}$x\qff \neq\qff y$\nsp.
{\off}Let\qss $D\qff =\qff (S\qff +\qff x)^*$\dnsp.\oss
Then\qss $\tau_{\fff x\fff y}\fff(D)\dff\in\dff\bbb^*$\dss
and exactly one of the elements\qss $x\comm y$ {\qff}is\dss in {\qff}$D$\dfdot}\newpage 

\proof\qss Since\qss $S\qff +\qff x\dff\in\dff\bbb$\qss 
and\qss $\tau_{\fff x\fff y}\fff(S\qff +\qff x)\qff =\qff S\qff +\qff y\dff\in\dff\bbb$\dnsp,\oss
we see that
\begin{equation}
\label{star-tau}
\quad
\tau_{\fff x\fff y}\fff(\dff(S\qff +\qff x)^*)\off 
=\off (\fff\tau_{\fff x\fff y}\fff(S\qff +\qff x)\dff)^*\off
=\off (S\qff +\qff y)^*\dff\in\dff\bbb^*.
\end{equation}
by the condition $\mathbold{L1}$\dfdot 
Obviously{\nsp},\qss $S\qff +\qff x\qff \neq\qff S\qff +\qff y$\dfdot
Therefore
\begin{equation}
\label{star-neq}
\quad
(S\qff +\qff y)^*\off
\neq  (S\qff +\qff x)^*\
\end{equation}
by the injectivity of the linking map.\qss
Since\qss $D\qff =\qff (S\qff +\qff x)^*$\dfcom
{\qff}(\ref{star-tau})\qff\oss implies that\qss $\tau_{\fff x\fff y}\fff(D)\dff\in\dff\bbb^*$\dfcom
and\qss (\ref{star-tau})\qss and\qss (\ref{star-neq})\qss together imply that\qss 
$\tau_{\fff x\fff y}\fff(D)\qff \neq\qff D$\dfdot
Finally{\nsp},\qss $\tau_{\fff x\fff y}\fff(D)\qff \neq\qff D$\qss implies that
exactly one of the elements\qss $x\comm y$ {\qff}is\dss in {\qff}$D$\dfdot  \eproof

\myitpar{Proof of Lemma \ref{xy-u-sets}.}\qss Let\qss $D\qff =\qff (S\qff +\qff x)^*\qff =\qff A\qff +\qff y$\dfdot
By Lemma\qss \ref{xy-flip}\qss $\tau_{\fff x\fff y}\fff(D)\dff\in\dff\bbb^*$\qss and exactly one of the elements
$x\comm y$ is in $D$\dfdot
Since\qss $y\dff\in\dff A\qff +\qff y\qff =\qff D$\dfcom
in fact\qss $x\dff\not\in\dff D$\qss and\qss $y\dff\in\dff D$\dnsp.\oss
It follows that 
$A\qff +\qff x\qff 
=\qff \tau_{\fff x\fff y}\fff(A\qff +\qff y)\qff 
=\qff \tau_{\fff x\fff y}\fff(D)\dff\in\dff\bbb^*$\dfcom
and hence $x\dff\in\dff U^*(A)$\dfdot  \eproof

\myappend{Classification\qss of\qss linkings}{app-linkings}

\vspace*{\bigskipamount}
Let\qss $\bbb$\dfcom $\bbb^*$\qss be pre-matroids on\qss $X$\ffcom
and let\dff\qss $L\dff\colon \bullet\qff\longmapsto\qff\bullet^*$\dff\qss 
be a linking\qss $\bbb\qff\longmapsto\qff\bbb^*$\dnsp.\qff\oss
Recall that by\dss $\tau_{\fff a b}$\dss we denote the transposition of distinct elements\dss 
$a\fff,\off b\dff\in\dff X$\dnsp.\oss

\myappar{Lemma.}{tau-fixed}\oss \emph{Suppose that\dff\qss $B\dff \in\dff \bbb$\qss
and\qss $x\fff,\off y$\qss are two different elements of\oss $X$\nsp.\qff\oss 
If\oss $x\dff \in\dff B^*$\qss and\qss $\tau_{\dff x\dff y}\fff(B)\qff =\qff B$\dnsp,\qff\oss
then\qss $y\dff \in\dff B^*$\dfdot}

\proof\qss By applying\qss $\bullet\qff \mapsto\qff \bullet^*$\qss to\qss 
$\tau_{\dff x\dff y}\fff(B)\off =\off B$\dnsp,\qff\oss 
we see that\qss
$\tau_{\dff x\dff y}\fff(B)^*\off =\off B^*$\dnsp.\qff\oss
Since\qss $\tau_{\dff x\dff y}\fff(B)\off =\off B\dff \in\dff \bbb$\dnsp,\qff\oss
the property\dss $\mathbold{L1}$\dss implies that\oss 
$\tau_{\dff x\dff y}\fff(B^*)\off 
=\off \tau_{\dff x\dff y}\fff(B)^*\off 
=\off B^*
$\dnsp.\qff\oss
Since\qss $x\dff \in\dff B^*$\oss and\oss $\tau_{\dff x\dff y}\fff(x)\qff =\qff y$\nsp,\qff\oss 
this implies\qss $y\dff \in\dff B^*$\dnsp.\qss  \eproof

\myappar{Lemma.}{id-or-comp-basis} \emph{For every\qss $B\in\bbb$\qss 
either\qss $B^*\qff =\qff B$\dnsp,\qff\oss or\dff\oss $B^*\qff =\qff B^{\ccomp}$\dfdot}

\proof\qss 
\mytitle{Case 1\dff:}\qff\oss $B\cap B^*\off \neq\off \varnothing$\dfdot\oss
Let us fix an element\qss $x\in B\cap B^*$\dfcom
and let\qss $y$\qss be an arbitrary element of\qss $B-x$\nsp.\qff\oss
Then\qss $\tau_{\dff x\dff y}\fff(B)\qff =\qff B$\dnsp,\qff\oss
and hence Lemma\qss \ref{tau-fixed}\qss implies that\qss $y\dff \in\dff B^*$\dnsp.\oss
It follows that\qss $B-x\qff \subset\qff B^*$\dnsp,\oss and hence\qss $B\qff \subset\qff B^*$\dnsp.\oss
By applying this argument to the inverse linking\qss $\bbb^*\qff\longrightarrow\qff \bbb$\dnsp,\oss
we see that\qss $B^*\qff \subset\qff B$\qss also,\oss 
and hence\qss $B^*\off =\off B$\dfdot
\newpage

\myitpar{Case 2\dff:}\qff\oss $B\cap B^*\off =\off \varnothing$\dfdot\quad
Then\qss $B^*\qff \subset\qff B\halfff\csup$\dnsp.\qff\oss
Let us fix an element\qss $x\in B^*$\dfcom
and let\qss $y$\qss be an arbitrary element of\qss $B\halfff\csup -x$\nsp.\qff\oss
The both $x\fff,\off y\dff \not\in\dff B$\dnsp,\oss
and hence\qss $\tau_{\dff x\dff y}\fff(B)\qff =\qff B$\dnsp.\qff\oss
Now Lemma\qss \ref{tau-fixed}\qss implies that\qss $y\dff \in\dff B^*$\dnsp.\oss
As in Case 1,\oss it follows that\qss 
$B\halfff\csup -x\qff \subset\qff B^*$\qss and hence\qss $B\halfff\csup\qff \subset\qff B^*$\dnsp.\oss
Since\qss $B^*\qff \subset\qff B\halfff\csup$\dnsp,\oss 
it follows that in this case\qss $B^*\off =\off B\halfff\csup$\dfdot   \eproof

\myappar{Lemma.}{basis-to-basis} \emph{Suppose that\qss $\tau$\qss is a transposition.\oss
If\oss $B^*\off =\off B$\dnsp,\oss then\qss $\tau\dff(B)^*\off =\off \tau\dff(B)$\dnsp,\qff\oss
and\dss if\oss  $B^*\off =\off B\fff\csup$\dnsp,\qff\oss 
then\oss $\tau\dff(B)^*\off =\off \tau\dff(B)\fff\csup$\dnsp.}

\proof\qss Note\dss that\qss $\tau\dff(B)^*\off =\off \tau\dff(B^*)$\qss
by the linking property\qss $\mathbold{L1}$\dnsp,\qff\oss and\oss
$\tau\dff(B\fff\csup)\off =\off \tau\dff(B)\fff\csup$\oss 
because\qss $\tau$\qss is a bijection.\qff\oss  
Therefore,\oss if\qss $B^*\off =\off B$\dnsp,\oss then\oss
\[
\quad
\tau\dff(B)^*\off =\off \tau\dff(B^*)\off =\off \tau\dff(B),
\]
and if\oss $B^*\off =\off B\fff\csup$\dnsp,\qff\oss then\oss 
$\displaystyle
\tau\dff(B)^*\off =\off \tau\dff(B^*)\off =\off \tau\dff(B\fff\csup)\off 
=\off \tau\dff(B)\fff\csup$\dfdot\oss  \eproof

\myappar{Lemma.}{chains-of-bases} \emph{Suppose that\qss $B\fff,\off B'$\qss are two bases of\oss $\bbb$\dnsp.\qff\oss
Then there exists a sequence
\[
\quad
B\qff =\qff B_{\dff 1}\fff,\hspace*{1.0em} B_{\dff 2}\fff,\hspace*{1.0em} \ldots\fff,\hspace*{1.0em} 
B_{\dff n}\qff =\qff B'
\]
of bases of\oss $\bbb$\oss such that\oss $B_{\dff i+1}\off =\off \tau_{\dff i}\dff(B_{\dff i})$\off
for some transposition\qss $\tau_{\dff i}$\qss
for every\qss $i = 1\fff,\off 2\fff,\off \ldots\fff,\off n-1$\ffdot}

\proof\qss If\qss $B'\off =\off B$\dfcom there is nothing to prove.\qff\oss
If\qss $B'\off \neq\off B$\dnsp,\qff\oss 
then the symmetric difference\qss $B\sdiff B'$\qss is not empty{\nsp}.\qss\oss
By\dss Theorem\qss \ref{symm-exchange},\oss 
in this case there are two elements $x\dff\in\dff B\qff \smallsetminus\qff B'$\dnsp,\qff\oss
$y\dff\in\dff B'\qff \smallsetminus\qff B$\qss such that\qss 
$B\sdiff\{\dff x\fff,\off y\dff\}\dff\in\dff\bbb$\dnsp.\qff\oss
Let 
\[
\quad
B_1\off =\off B\sdiff\{\dff x\comm y\dff\}\off 
=\off B\qff -\qff x\qff +\qff y\fff.
\]
Then\oss $B_1\dff \cap\dff B'\off =\off (\fff B\dff \cap\dff B'\fff)\qff +\qff y$\oss
because\oss $B\sdiff\{\dff x\fff,\off y\dff\}\off =\off B\qff -\qff x\qff +\qff y$\dnsp.\oss 
It follows that the intersection\oss $B_1\dff \cap\dff B'$\oss 
contains one element more than\oss $B\dff \cap\dff B'$\dnsp.\qff\oss
An application of the induction completes the proof.\qss   \eproof

\myitpar{Proof of Theorem \ref{linkings-classification}.}\qss 
Let $B_1\in\bbb$\dfdot 
By Lemma\qss \ref{id-or-comp-basis}\qss for every basis\qss $B$\qss of\qss $\bbb$\qss
either\qss $B^*\qff =\qff B$\dnsp,\oss or\oss $B^*\qff =\qff B\fff\csup$\dnsp.\oss

Lemma\oss \ref{basis-to-basis}\oss together with Lemma\oss \ref{chains-of-bases}\oss imply 
that if\oss
$B^*\off =\off B$\oss 
for one basis\qss 
$B\dff\in\dff\bbb$\dnsp,\hspace*{0.7em}
then\hspace*{0.5em} 
$B^*\off =\off B$\dff\oss 
for all bases\dff\oss $
B\dff\in\dff\bbb$\dnsp,\qff\oss
and\dss if\oss
$B^*\off =\off B\fff\csup$\oss 
for one basis\qss 
$B\dff\qff\in\qff\bbb$\dnsp,\hspace*{0.85em}
then\hspace*{0.6em} 
$B^*\dff\off =\qff\off B\fff\csup$\oss 
for all bases\qff\oss 
$B\dff\in\dff\bbb$\dnsp.\qss  \eproof

\mynonumsection{Note\qss historique}

\vspace*{\bigskipamount}
W.T. Tutte was uncommonly generous in sharing
the route which lead him to his discoveries.\qss
He was equally generous in giving the credit to his predecessors.\qss
In particular{},\qss he described the way which lead him to the definition of
the Tutte polynomial and to the Tutte order-independence theorem in Chapter 5 of his book\qss \cite{t3},\qss
and,\qss with additional details and from a somewhat different prespective,\qss
in a remarkable paper\qss \cite{t4}.\qss

Together with his coauthors on the paper\qss \cite{bsst},\qss
Tutte observed that the number $C\dff(G)$ of the spanning trees (or forests) 
of a graph $G$ satisfies the relation
\begin{equation}
\label{del-con}
\quad
C\dff(G)\off =\off C\dff(G'_A)\qff +\qff C\dff(G''_A),
\end{equation}
where $A$ is an arbitrary edge of $G$\dfcom
$G'_A$ is the result of the deleting $A$ from $G$\dfcom
and $G''_A$ is the result of contracting $A$ together with its two endpoints into a single vertex.\qss
Tutte discovered his polynomial while looking for invariants of graphs 
satisfying relations similar to (\ref{del-con}).\qss
As Tutte wrote in\qss \cite{t4},\qss he\qss 
\emph{``...come across one such in a footnote to one of Hassler Whitney's papers''}\qss \cite{w2}.\qss

This footnote is,\qss in fact,\qss the note added in proofs at the very end of\qss \cite{w2}.\qss
Whitney\qss \cite{w1} introduced an invariant $m_{\dff i\dff j}$ 
of graphs for every pair of non-negative integers $i\comm j$\dfdot
The paper\qss \cite{w1} is devoted to the proof of what is now known as the inclusion-exclusion principle,\qss
but was called by Whitney\qss \emph{the logical expansion},\qss
and to various applications of it.\qss
One of the application is to the number of colorings of a graph.\qss
The inclusion-exclusion principle applied to colorings inevitably leads 
to the Whiteny invariants $m_{\dff i\dff j}$\dfdot
In \cite{w2} Whitney continued to study these invariants and their applications to the colorings of graphs.\qss
These invariants satisfy the relation
\begin{equation}
\label{del-con-mij}
\quad
m_{\dff i\dff j}\dff(G)\off =\off m_{\dff i\dff j}\dff(G'_A)\qff +\qff m_{\dff i-1,\dff j}\dff(G''_A).
\end{equation}
Apparently{\nsp},\qss this is an observation of\qss R.M. Foster{},\qss 
who used it to calculate $m_{\dff i\dff j}$ 
for large number of graphs,\qss according to\qss \cite{w2}.\qss 

In \cite{t1},\qss Tutte developed a general theory of graph invariants statisfying recursion relations 
similar to\qss (\ref{del-con}),\qss (\ref{del-con-mij}).\qss 
He called them $W$\dnsp-\emph{functions}.\qss
Let us now quote\qss \cite{t4}.\qss 
\vspace*{-\bigskipamount}\vspace*{\medskipamount}
\begin{quote}
Playing with my $W$\dnsp-functions I obtained a two-variable polynomial...

... In my papers I called this function the dichromate,\qss
but it is now generally known as the Tutte polynomial.\qss
This may be unfair to Hassler Whitney who knew and 
used analogous coefficients without bothering to affix 
them to two variables\qss \cite{w2}.
\end{quote}
\vspace*{-\bigskipamount}\vspace*{\medskipamount}
From the point of view adopted in the present paper{},\qss
there is no reason to attach the coefficients to a polynomial.\qss
Since not even the addition of polynomials is used,\qss
it is more natural to deal directly with the coefficients.\qss
It is convenient to arrange these coefficients (which are non-negative integers)
into a single entity{\nsp},\qss namely{\nsp},\qss into a multi-set.\qss
Multi-sets are a special case of\qss
\emph{generalized sets},\qss introduced by Whitney in the paper\qss \cite{w3},\qss
closely related to\qss \cite{w1}.\qss 

Tutte observed that the value of his polynomial $\chi\fff(G;\qff\mathbold{x},\dff\mathbold{y})$ 
at $(\mathbold{x},\dff\mathbold{y})\qff =\qff (1,\dff 1)$
is equal to the number of spanning trees of $G$\dfdot 
This lead him to the hypothesis that $\chi\fff(G;\qff\mathbold{x},\dff\mathbold{y})$
can be presented as a sum of\qss \emph{``something simple''}\qss over all spanning trees of $G$\dfdot
Initially{\nsp},\qss this seemed to be impossible because
even for very simple graphs $G$
the symmetries of $G$ are not reflected in\qss
$\chi\fff(G;\qff\mathbold{x},\dff\mathbold{y})$\dfdot
The way out,\qss found by Tutte,\qss 
was to break any potential symmetry by enumerating the edges of $G$\dfdot
This lead him to his remarkable definition of  
$\chi\fff(G;\qff\mathbold{x},\dff\mathbold{y})$ 
in terms of the internal and external activities\qss
(see Section \ref{orders-pre-matroids}),\qss
and to his order-invariance theorem.\qss
Let us quote Tutte\qss \cite{t3}\qss again.
\vspace*{-\bigskipamount}\vspace*{\medskipamount}
\begin{quote}
I marvelled that all the different possible enumerations 
should give rise to the same polynomial 
$\chi\fff(G;\qff\mathbold{x},\dff\mathbold{y})$\dfcom
even though different enumerations usually gave 
different internal and external activities for a given spanning tree.\qss 
But I recalled that Hassler Whitney,\qss 
giving the chromatic polynomial in terms of broken circuits,\qss 
had encountered a similar phenomenon\qss \cite{w1}.
\end{quote}
\vspace*{-\bigskipamount}\vspace*{\medskipamount}
It seems that the role of contributions of Hassler Whitney 
to the modern combinatorics in general,\qss
and to theory of polynomial invariants of graphs and matroids in particular{},\qss 
is at the very least under-appreciated.\qss
Whitney also is one of discoverers of matroids\qss \cite{w4},\qss 
which provide the proper context for the latter theory.\qss
On the other hand,\qss the beauty and originality of 
the definition of the Tutte polynomial in terms of the enumerations of edges,\qss
as well as the Tutte order-independence theorem,\qss are also under-appreciated.\qss

The relations of the form\qss (\ref{del-con}),\qss (\ref{del-con-mij}),\qss
apparently{\nsp},\qss dominated the field from the very beginning.\qss
This domination was later reinforced by a somewhat superficial analogy with 
the Grothendieck construction of $K$\dnsp-groups in terms of generators and relations.\qss
In fact,\qss similar constructions were used before Grothendieck.\qss 
The easiest example is the standard construction of 
the integers $\zzz$ from the natural numbers $\nnn$\dfdot
More importantly{},\qss in the framework of this analogy 
graphs and matroids correspond to vector bundles on 
a\qss \emph{fixed}\qss algebraic variety{\nsp},\qss
and the abelian group $\zzz\dff[\dff\mathbold{x}\dff,\qff\mathbold{y}\dff]$
corresponds to the $K$\dnsp-group of this variety{\nsp}.\qss
Grothendieck $K$\dnsp-groups form a functor on the category of algebraic varieties,\qss
but there is no corresponding category in the theory of the Tutte polynomials.\qss

In more recent times,\qss relations similar to\qss (\ref{del-con}),\qss (\ref{del-con-mij}) 
come to prominence in topology and led to new invariants of knots and $3$\dnsp-dimensional manifolds.\qss
The analogy with these invariants is much closer than with the Grothendieck $K$\dnsp-groups.\qss
As in the theory of the Tutte polynomial,\qss these invariants  
do not lead to a functor similar to the Grothendieck $K$\dnsp-functor{}.\qss 
Moreover{},\qss these invariants turned out to be 
directly related with the Tutte polynomial.\qss
See,\qss for example,\qss M.B. Thistlethwaite's paper \cite{th}.\qss
Unfortunately{\nsp}, these remarkable developments moved 
the ideas of Tutte \cite{t2} even further into the shadows.\qss

\vspace*{3ex}

\begin{flushright}

April 12,\qss 2016

\vspace{\bigskipamount}
 
http:/\!/\hspace*{-0.07em}nikolaivivanov.com\\\vspace{\bigskipamount}
E-mail\fff:\qff\qss nikolai@nikolaivivanov.com

\end{flushright}

\end{document}